\documentclass[11pt,a4paper]{article}
\usepackage[T1]{fontenc}
\IfFileExists{lmodern.sty}{\usepackage{lmodern}}{}
\usepackage[margin=27mm,headheight=15pt]{geometry}
\usepackage{amsmath,amssymb,amsthm}
\usepackage{fancyhdr,needspace,tikz-cd}
\usepackage[pdfborder={0 0 0},hyperfootnotes=false]{hyperref}
\DeclareRobustCommand{\modu}[1]{\ensuremath{#1\text{-}\mathsf{mod}}}
\DeclareRobustCommand{\opp}{^{\mathrm{op}}}
\DeclareMathOperator{\Hom}{Hom}
\DeclareMathOperator{\End}{End}
\DeclareMathOperator{\Ext}{Ext}
\DeclareMathOperator{\Tr}{Tr}
\DeclareMathOperator{\Image}{Im}
\DeclareMathOperator{\Ker}{Ker}
\DeclareMathOperator{\Coker}{Coker}
\DeclareMathOperator{\soc}{soc}
\DeclareMathOperator{\Head}{top}
\DeclareMathOperator{\ann}{ann}
\DeclareMathOperator{\add}{add}
\DeclareMathOperator{\Fac}{Fac}
\DeclareMathOperator{\Sub}{Sub}
\DeclareMathOperator{\gldim}{gl.dim}
\DeclareMathOperator{\repdim}{repdim}

\newtheorem{lemma}{Lemma}[section]
\newtheorem{proposition}[lemma]{Proposition}
\newtheorem{corollary}[lemma]{Corollary}
\newtheorem{theorem}[lemma]{Theorem}
\theoremstyle{remark}
\newtheorem{remark}[lemma]{Remark}
\theoremstyle{plain}
\numberwithin{equation}{section}

\title{\bfseries Auslander-Reiten translation and the representation dimension of tame algebras}
\author{Liang Chen \quad Zhenxian Chen \quad Yongyi Wen$^{*}$ \quad Weiheng Xia}
\date{}
\hypersetup{pdftitle={Auslander-Reiten translation and the representation dimension of tame algebras},pdfauthor={Liang Chen, Zhenxian Chen, Yongyi Wen, Weiheng Xia},pdfkeywords={Tame representation type, Auslander-Reiten translate, Representation embedding, Representation dimension, Homological conjectures}}
\begin{document}
\maketitle
\begingroup
\renewcommand{\thefootnote}{}
\footnotetext{\textsuperscript{*}\,Corresponding author. Email:
\href{mailto:2240501002@cnu.edu.cn}{2240501002@cnu.edu.cn}.
\par\noindent\hspace*{1.8em}2020 Mathematics Subject Classification: Primary 16G60; Secondary 16G70, 16E10.
\par\noindent\hspace*{1.8em}Keywords: Tame representation type; Auslander--Reiten translate;
Representation embedding; Representation dimension; Homological conjectures.}
\endgroup

\begin{abstract}
We prove that a finite-dimensional algebra is tame if and only if, in each dimension, Auslander--Reiten
translation fixes all but finitely many isomorphism classes of
indecomposable modules. This also proves Crawley-Boevey's conjecture that
homogeneous tubes characterize tame representation type. We also prove
that every tame algebra has representation dimension at most three.
Consequently, tame algebras satisfy the Auslander--Reiten conjecture
and the finitistic dimension conjecture.
\end{abstract}

\tableofcontents

\section{Introduction}

How can one recognize tame representation type from the
Auslander--Reiten quiver? Crawley-Boevey \cite[Theorem~D and Corollary~E]{CB} gave two
necessary conditions and conjectured that each is sufficient
\cite[p.~483]{CB}. We prove both converses and we also show that tame algebras have
representation dimension at most three.

Unless otherwise stated, $k$ is an algebraically closed field and $A$
is a finite-dimensional $k$-algebra. All algebras are associative with identity. Unless a side is specified,
modules are left modules. For a $k$-algebra $R$, we write $\modu R$
for the category of finite-dimensional $R$-modules.

An additive category is of \emph{finite representation type} if it has
only finitely many isomorphism classes of indecomposable objects. We
apply the same terminology to $A$ when $\modu A$ has this property.
A \emph{representation embedding} between $k$-linear abelian
Krull--Schmidt categories is an exact $k$-linear functor $E$ that preserves
indecomposable objects and reflects isomorphism classes. The latter
condition means that $E(X)\cong E(Y)$ implies $X\cong Y$.

We recall the standard definitions of wild and tame representation type
(see \cite{Drozd1977} and \cite[Definitions~6.4 and~6.5]{CB}).
Put $S=k\langle x,y\rangle$, the free algebra on two generators.
The algebra $A$ is \emph{of wild representation type}, or \emph{wild},
if there is an $A$-$S$-bimodule $T$, free of positive finite rank as
a right $S$-module, for which
$T\otimes_S-\colon\modu S\to\modu A$ is a representation embedding.
The algebra $A$ is \emph{of tame representation type}, or \emph{tame},
if for each $d\geq1$ there are finitely many $A$-$k[t]$-bimodules
$Q_1,\ldots,Q_m$, free of rank $d$ on the right, with the following
property. Up to isomorphism, all but finitely many indecomposable
$A$-modules of dimension $d$ have the form
$Q_j\otimes_{k[t]}k[t]/(t-a)$ for some $j$ and $a\in k$. Each exceptional module $X$ can be included
by adding the constant family $X\otimes_k k[t]$. The list may be empty,
so this convention includes algebras of finite representation type.

An indecomposable $A$-module $X$ is \emph{homogeneous} if
$X\cong\tau_A(X)$, where $\tau_A$ is the Auslander--Reiten translate.
We call $A$ \emph{$\tau$-tame} if, for each $d\geq1$, there are only
finitely many isomorphism classes of $d$-dimensional indecomposable
$A$-modules $X$ with $X\not\cong\tau_A(X)$. By Crawley-Boevey's theorem
\cite[Theorem~D]{CB}, every tame algebra is $\tau$-tame.
Nagase \cite[Theorem~4.1]{NagaseLocal} treated wild local algebras
and certain two-point algebras. See also \cite{Nagase}. Zhang and Xu stated the converse of Theorem~D in
\cite[Theorem~3]{ZX}.
Our first theorem gives a shorter proof using self-extensions and
bounds the least dimension in
which a wild algebra has infinitely many pairwise nonisomorphic
nonhomogeneous indecomposable modules.

\begin{theorem}\label{thm:main}
Let $T$ be an $A$-$S$-bimodule, free of rank $r\geq1$ as a right
$S$-module. Suppose that $T\otimes_S-$ is a representation embedding.
Then for some integer $d$ with $1\leq d\leq r^2$ there are infinitely
many pairwise nonisomorphic indecomposable $A$-modules $X$ of dimension
$d$ with $X\not\cong\tau_A(X)$. In particular, $A$ is tame if and only
if $A$ is $\tau$-tame.
\end{theorem}

A homogeneous tube is a stable tube of rank one in the
Auslander--Reiten quiver. The following corollary proves the converse
of Crawley-Boevey's Corollary~E, which de la Pe\~na also listed as
an open problem \cite[Section~8.3.8, p.~207]{deLaPena}.

\begin{theorem}\label{thm:main2}
The algebra $A$ is tame if and only if, for each $d\geq 1$, all but
finitely many isomorphism classes of $d$-dimensional indecomposable
$A$-modules lie in homogeneous tubes.
\end{theorem}

The second question concerns the representation dimension of a tame
algebra. Recall that $\repdim A$ is the minimum of $\gldim\End_A(G)$ over all
generator-cogenerators $G$ in $\modu A$. A module $G$ is a
generator-cogenerator if every indecomposable projective or injective
$A$-module is isomorphic to a direct summand of $G$.
Auslander \cite{Auslander} proved that $A$ is of finite representation
type if and only if $\repdim A\leq2$. The conjecture that
$\repdim A\leq3$ for every tame algebra was recorded by Oppermann
\cite[Conjecture~1.3]{Oppermann}. Erdmann, Holm, Iyama and Schr\"oer
\cite[Corollary~1.3]{EHIS} proved this bound for special biserial
algebras. Assem, Skowro\'nski and Trepode \cite{AST} proved that
self-injective algebras of Euclidean type have representation dimension
three. Theorem~\ref{thm:B} proves this conjecture over algebraically closed fields.

\begin{theorem}\label{thm:B}
If $A$ is of tame representation type, then $\repdim A\leq3$.
\end{theorem}

Theorem~\ref{thm:B} has several homological consequences.
The little finitistic dimension conjecture asks whether the finite
projective dimensions of finitely generated modules have a uniform
upper bound. Corollary~\ref{cor:findim} establishes this conjecture for every
tame algebra and its opposite. The Auslander--Reiten conjecture
\cite{AR1975} asserts that a module $M\in\modu A$ is projective whenever
$\Ext_A^i(M,M\oplus A)=0$ for all $i>0$. The Gorenstein symmetry
conjecture asserts that the regular module has finite injective
dimension on the left if and only if it has finite injective dimension
on the right. Corollaries~\ref{cor:projectivity}
and~\ref{cor:gorenstein-symmetry} establish both conjectures for tame
algebras.

We write morphisms on the right and functors on the left. The composite
of $f\colon X\to Y$ and $g\colon Y\to Z$ is $fg$, with
$(x)(fg)=((x)f)g$. The composite of functors $F$ followed by $G$ is
$GF$. We write $J(R)$ for the Jacobson radical of a ring $R$,
$\overline{x}$ for the image of $x$ in a quotient, and
$D=\Hom_k(-,k)$ for the standard duality between $\modu A$ and
$\modu{A\opp}$. For background on almost split sequences we refer to
\cite{ARS}.

The paper is organized as follows. Section~\ref{sec:extension-rank}
uses Auslander--Reiten duality and change of coefficients to bound the
rank of a map into a group of self-extensions. Section~\ref{sec:wild}
constructs independent self-extensions over the free algebra and uses
their images to obtain a lower bound for the same rank. Comparing the
two bounds proves Theorem~\ref{thm:main}, from which
Theorem~\ref{thm:main2} follows.
Section~\ref{sec:ext} compares extensions over $A$ and over a quotient
by a power of its radical. In Section~\ref{sec:dimensions}, the
resulting bound and Lemma~\ref{lem:last-radical-layer} give finite
representation type for the category of quotients of finite
sums of the radical. We then construct the right approximations for
Theorem~\ref{thm:B} and deduce its homological applications.

\section{Ranks of extension maps}\label{sec:extension-rank}

Crawley-Boevey's homogeneity theorem \cite[Theorem~D]{CB} shows
that, for a tame algebra, all but finitely many indecomposable modules
in each dimension are isomorphic to their Auslander--Reiten translates.
For Theorem~\ref{thm:main}, we use the corresponding finiteness
condition only in dimensions at most $d_0$. Under this finiteness condition, Proposition~\ref{prop:upper} bounds
the rank of the extension map defined in Subsection~\ref{sec:rank}.

\subsection{The extension map}\label{sec:rank}

\noindent\textbf{In this subsection,} let $k$ be an arbitrary field
and let $\Lambda$ be a $k$-algebra. We use
Yoneda extensions. For $\Lambda$-module homomorphisms $a\colon X'\to X$
and $b\colon Y\to Y'$,
$\Ext_\Lambda^1(a,Y)$ sends an extension class to the class of the pullback
along $a$, and $\Ext_\Lambda^1(X,b)$ sends an extension class to the class
of the pushout along $b$.

Let $V$ be a finite-dimensional $k$-module. Define
$C=k\oplus V$ by $(a,v)(a',v')=(aa',av'+a'v)$. Then $V^2=0$,
and $\kappa=C/V$ is the unique simple $C$-module up to isomorphism.
We first describe the self-extensions of $\kappa$ in terms of
$V^*=\Hom_k(V,k)$. For $\ell\in V^*$, the map
$v\mapsto ((v)\ell)\overline{1}$ from $V$ to $\kappa$ is a
$C$-module homomorphism, since $V$ annihilates both modules.
Pushing out $0\to V\to C\to\kappa\to0$ along this map gives
\[
 \xi_\ell\colon\quad 0\longrightarrow\kappa\xrightarrow{j_\ell}
 N_\ell\xrightarrow{p_\ell}\kappa\longrightarrow0.
\]
Here $N_\ell$ is the quotient of $\kappa\oplus C$ by the submodule
formed by the pairs $(((v)\ell)\overline{1},-v)$ with $v\in V$.
Let $n_1$ and $n_2$ be the classes of $(\overline{1},0)$ and $(0,1)$,
respectively. The elements $n_1,n_2$ form a $k$-basis of $N_\ell$, and
\[
 (\overline{1})j_\ell=n_1,\qquad (n_2)p_\ell=\overline{1},
 \qquad vn_1=0,\qquad vn_2=((v)\ell)n_1\quad(v\in V).
\]
Applying $\Hom_C(-,\kappa)$ to $0\to V\to C\to\kappa\to0$ gives
the exact sequence
\[
 \Hom_C(C,\kappa)\longrightarrow\Hom_C(V,\kappa)
 \longrightarrow\Ext^1_C(\kappa,\kappa)\longrightarrow0.
\]
For every $f\colon C\to\kappa$ and $v\in V$ we have
$(v)f=v((1)f)=0$. Thus the first homomorphism is zero. Composing
$\ell\mapsto(v\mapsto ((v)\ell)\overline{1})$ with the connecting
homomorphism gives the isomorphism
$V^*\xrightarrow{\sim}\Ext^1_C(\kappa,\kappa)$,
$\ell\mapsto[\xi_\ell]$.

We consider $\Lambda$-$C$-bimodules that are free of finite rank as
right $C$-modules.
For each such $X$, put $X_0=X\otimes_C\kappa$. Since
$X\otimes_C-$ is exact, we may define
\[
 d_X\colon V^*\longrightarrow\Ext^1_\Lambda(X_0,X_0),\qquad
 \ell\longmapsto [X\otimes_C\xi_\ell].
\]
The functor $X\otimes_C-$ preserves the direct sums, pullbacks and
pushouts used to form Baer sums. Since $1_X\otimes_C(a1_\kappa)=a1_{X_0}$
for $a\in k$, tensoring also preserves scalar multiplication on extension
classes. Thus $d_X$ is a homomorphism of $k$-modules. Write
$\Delta(X)=\dim_k\Image(d_X)$. In particular, $\Delta(X)$ is defined
when $\Lambda=k\langle x,y\rangle$.

For a bimodule homomorphism $a\colon X\to Y$, put
$a_0=a\otimes_C 1_\kappa\colon X_0\to Y_0$. For each $C$-module $N$,
the component of $a\otimes_C-\colon X\otimes_C-\to Y\otimes_C-$ at $N$
sends $x\otimes n$ to $(x)a\otimes n$. Naturality for $j_\ell$ and
$p_\ell$ gives the following diagram of $\Lambda$-modules with exact rows.
All tensor products in the diagram are over $C$.
\[
\begin{tikzcd}[column sep=1.9em,row sep=2.1em]
0 \arrow{r}
 & X_0 \arrow{r}{1_X\otimes j_\ell} \arrow{d}[swap]{a_0}
 & X\otimes_C N_\ell \arrow{r}{1_X\otimes p_\ell}
     \arrow{d}{a\otimes 1_{N_\ell}}
 & X_0 \arrow{r} \arrow{d}{a_0} & 0 \\
0 \arrow{r}
 & Y_0 \arrow{r}[swap]{1_Y\otimes j_\ell}
 & Y\otimes_C N_\ell \arrow{r}[swap]{1_Y\otimes p_\ell}
 & Y_0 \arrow{r} & 0.
\end{tikzcd}
\]
Pushing out the top row along $a_0$ and pulling back the bottom row
along $a_0$ give the same extension class. Hence
\begin{equation}\label{eq:rank-naturality}
 \bigl((\ell)d_X\bigr)\Ext_\Lambda^1(X_0,a_0)
 =\bigl((\ell)d_Y\bigr)\Ext_\Lambda^1(a_0,Y_0)
 \qquad\text{in }\Ext_\Lambda^1(X_0,Y_0).
\end{equation}
If $a$ is an isomorphism, then so is $a_0$. Both homomorphisms
on extension groups in \eqref{eq:rank-naturality} are then
isomorphisms, and $\Delta(X)=\Delta(Y)$.

For a finite-dimensional $\Lambda$-module $U$, the natural isomorphism
$(U\otimes_k C)\otimes_C N\to U\otimes_k N$ sends $(u\otimes c)\otimes n$
to $u\otimes cn$. Choose a $k$-module section $s\colon\kappa\to N_\ell$
of $p_\ell$. Then $1_U\otimes_k s$ is a $\Lambda$-module section of
$1_U\otimes_k p_\ell$. Thus $\Delta(U\otimes_k C)=0$.

Tensoring $\xi_\ell$ with $X\oplus Y$ gives the direct sum of
$X\otimes_C\xi_\ell$ and $Y\otimes_C\xi_\ell$. The two off-diagonal
components of this extension class are zero. Projection onto the
diagonal components therefore gives a monomorphism
\[
 \Image(d_{X\oplus Y})\hookrightarrow\Image(d_X)\oplus\Image(d_Y),\qquad
 (\ell)d_{X\oplus Y}\longmapsto\bigl((\ell)d_X,(\ell)d_Y\bigr).
\]
Hence $\Delta(X\oplus Y)\leq\Delta(X)+\Delta(Y)$.

For the rest of this subsection, assume that $\Lambda$ is
finite-dimensional and work in $\modu{\Lambda}$.
Let $P_1\xrightarrow{u}P_0\to X\to0$ be a minimal projective
presentation of a $\Lambda$-module $X$. The transpose
$\Tr_\Lambda(X)$ is the cokernel of
$\Hom_\Lambda(u,\Lambda)\colon\Hom_\Lambda(P_0,\Lambda)
\to\Hom_\Lambda(P_1,\Lambda)$. Write
$\tau_\Lambda(X)=D(\Tr_\Lambda(X))$, or simply $\tau(X)$ when the algebra
is clear. The Auslander--Reiten translate gives an equivalence from the module
category modulo projectives to the module category modulo injectives.
In particular, $\tau_\Lambda(P)=0$ for each projective module $P$.
For an indecomposable nonprojective module $X$, the almost split sequence
ending in $X$ starts at $\tau_\Lambda(X)$ \cite[Chapters~IV--V]{ARS}.

For $\Lambda$-modules $X,Y$, the Auslander--Reiten formula gives an epimorphism
of $k$-modules, natural in $Y$,
\begin{equation}\label{eq:rank-ar}
 \lambda_{X,Y}\colon\Hom_\Lambda(Y,\tau(X))\longrightarrow D(\Ext_\Lambda^1(X,Y)).
\end{equation}
To obtain $\lambda_{X,Y}$, take the quotient of
$\Hom_\Lambda(Y,\tau(X))$ by the homomorphisms that factor through
injective modules and apply the Auslander--Reiten isomorphism
\cite{ARIII}. See also
\cite[Chapter~XI, Section~5, pp.~392--393]{ARS}. For $X=X'\oplus P$
with $P$ projective, additivity extends the formula for $X'$ to $X$,
since $\tau(P)=0$ and $\Ext_\Lambda^1(P,Y)=0$.
For $\zeta\in\Ext_\Lambda^1(X,Y)$ and $f\colon Y\to\tau(X)$, write
$\langle\zeta,f\rangle_{X,Y}=(\zeta)((f)\lambda_{X,Y})$.
If $u\colon Y\to Y'$ and $f'\colon Y'\to\tau(X)$, then naturality gives
\begin{equation}\label{eq:rank-ar-natural}
 \big\langle(\zeta)\Ext_\Lambda^1(X,u),f'\big\rangle_{X,Y'}
 =\langle\zeta,uf'\rangle_{X,Y}.
\end{equation}

The next lemma compares $\Delta(X)$ and $\Delta(W)$ by factoring
endomorphisms of $X_0$ through $W_0$. The bound depends on the
codimension of the subspace spanned by the resulting composites.

\begin{lemma}\label{lem:rank}
Let $X$ and $W$ be $\Lambda$-$C$-bimodules as above and suppose
$X_0\cong\tau(X_0)$. Let $J'$ be a $k$-submodule of $\End_\Lambda(X_0)$
of codimension $q$. Suppose $J'$ has a basis $h_1,\ldots,h_t$ with
the following property. For each $i$, there are bimodule homomorphisms
$a_i\colon X\to W$ and $b_i\colon W\to X$ such that
\[
 h_i=a_{i,0}b_{i,0},\qquad
 a_{i,0}=a_i\otimes_C1_\kappa\colon X_0\to W_0,\qquad
 b_{i,0}=b_i\otimes_C1_\kappa\colon W_0\to X_0.
\]
Then $\Delta(X)\leq\Delta(W)+q$.
\end{lemma}

\begin{proof}
Put $L=(\Ker(d_W))d_X$. Let $\ell\in\Ker(d_W)$ and write $\zeta=(\ell)d_X$.
By \eqref{eq:rank-naturality} applied to $a_i$, we have
\[
 (\zeta)\Ext_\Lambda^1(X_0,a_{i,0})
 =\bigl((\ell)d_W\bigr)\Ext_\Lambda^1(a_{i,0},W_0)=0.
\]
Pushing out along $b_{i,0}$ gives $(\zeta)\Ext_\Lambda^1(X_0,h_i)=0$.
The map $h\mapsto(\zeta)\Ext_\Lambda^1(X_0,h)$ is a $k$-module
homomorphism. Since the $h_i$ span $J'$, we obtain
$(\zeta)\Ext_\Lambda^1(X_0,h)=0$ for every $h\in J'$.

Choose an isomorphism $\theta\colon X_0\to\tau(X_0)$. Composition
with $\theta$ gives an isomorphism
\[
 \End_\Lambda(X_0)\longrightarrow\Hom_\Lambda(X_0,\tau(X_0)),\qquad
 h\longmapsto h\theta.
\]
Compose this isomorphism with $\lambda_{X_0,X_0}$ in \eqref{eq:rank-ar}
and the restriction of functionals to $L\subseteq\Ext_\Lambda^1(X_0,X_0)$.
The restriction is surjective by exactness of $D$, so the composite
is surjective. For
$h\in J'$ and $\zeta\in L$, equation~\eqref{eq:rank-ar-natural} gives
\[
 \langle\zeta,h\theta\rangle_{X_0,X_0}
 =\big\langle(\zeta)\Ext_\Lambda^1(X_0,h),\theta\big\rangle_{X_0,X_0}=0.
\]
We therefore have an epimorphism
\[
 \End_\Lambda(X_0)/J'\twoheadrightarrow D(L),\qquad
 h+J'\longmapsto\bigl(\zeta\mapsto\langle\zeta,h\theta\rangle_{X_0,X_0}\bigr).
\]
In particular, $\dim_k L\leq q$.

There is also an epimorphism
\[
 \Image(d_W)\longrightarrow\Image(d_X)/L,\qquad
 [W\otimes_C\xi_\ell]\longmapsto [X\otimes_C\xi_\ell]+L.
\]
Indeed, if $(\ell)d_W=(\ell')d_W$, then $\ell-\ell'\in\Ker(d_W)$,
so $(\ell-\ell')d_X\in L$. Thus the map is well defined and is
surjective by the definition of $\Image(d_X)$. Hence
$\Delta(X)-\dim_k L\leq\Delta(W)$, as required.
\end{proof}

\subsection{Change of coefficients}\label{sec:ring}

To apply Lemma~\ref{lem:rank}, we construct modules over $C$ from
modules whose coefficients lie in a field extension of $k$.
We return to the standing assumptions on $k$ and $A$.
For Lemma~\ref{lem:ring} and Proposition~\ref{prop:upper}, fix a
$k$-module $V$ with basis $v_1,\ldots,v_m$, where $m\geq1$. As in Subsection~\ref{sec:rank},
put $C=k\oplus V$ and $\kappa=C/V$. Write $\rho\colon C\to k$ for
the scalar projection. Then $c+V\mapsto(c)\rho$ is a $k$-algebra
isomorphism $\kappa\to k$. Let $L$ be a separable closure of
$F=k(t_1,\ldots,t_m)$. The first step is to map a ring of
coefficients in $L$ to $C$ so that the image of each $t_i$ has
$V$-component $v_i$.

\begin{lemma}\label{lem:ring}
For finite subsets $\mathcal A\subseteq L$ and $\mathcal B\subseteq L^\times$,
there exist a finitely generated $k$-subalgebra $R\subseteq L$ and a
$k$-algebra homomorphism $\sigma\colon R\to C$ with the following
properties. The algebra $R$ contains $t_1,\ldots,t_m$ and
$\mathcal A\cup\mathcal B$, and $\sigma$ satisfies
\[
  (t_i)\sigma=a_i+v_i\quad(1\leq i\leq m),\qquad
  (b)\sigma\in C^\times\quad(b\in\mathcal B)
\]
for some $a_1,\ldots,a_m\in k$.
\end{lemma}

\begin{proof}
The elements of $\mathcal A\cup\mathcal B$, together with the inverses
of the elements of $\mathcal B$, generate a finite separable field
extension of $F$ inside $L$. By the primitive element theorem this
extension is $F(z)$ for some $z\in L$. Let $p\in F[Z]$ be the monic
minimal polynomial of $z$, and put $p'=\mathrm{d}p/\mathrm{d}Z$.
Since $p$ is separable, there exist $u,w\in F[Z]$ such that $up+wp'=1$.

Since $F(z)=F[z]$, for each
$x\in\mathcal A\cup\mathcal B\cup\{b^{-1}\mid b\in\mathcal B\}$ choose
$f_x\in F[Z]$ with $x=f_x(z)$. Choose a common denominator
$0\ne g\in k[t_1,\ldots,t_m]$ for the coefficients of $p,u,w$ and
all $f_x$. Thus
\[
 gp,\;gu,\;gw,\;gf_x\in k[t_1,\ldots,t_m][Z]
 \qquad\text{for every such }x.
\]
Put $R_0=k[t_1,\ldots,t_m,1/g]$ and $R=R_0[z]$. Then $R$ is a finitely
generated $k$-algebra, and $f_x\in R_0[Z]$ gives $x\in R$ for each
$x$ above. Evaluation gives the isomorphism
\[
 R_0[Z]/(p)\xrightarrow{\sim}R,\qquad h+(p)\longmapsto h(z).
\]
Indeed, division by the monic polynomial $p$ gives a remainder in
$R_0[Z]$ of degree less than $\deg p$. If $h(z)=0$, this remainder
has $z$ as a root and hence is zero by the minimality of $p$ over $F$.

Since $k$ is infinite, choose $a=(a_1,\ldots,a_m)\in k^m$ with $g(a)\ne0$.
Write $a+v=(a_1+v_1,\ldots,a_m+v_m)$. Evaluation at $a+v$ is a
$k$-algebra homomorphism $k[t_1,\ldots,t_m]\to C$. The difference
$g(a+v)-g(a)$ belongs to $V$. Since $V^2=0$, the inverse of $g(a+v)$ is
\[
 g(a)^{-1}-g(a)^{-2}\bigl(g(a+v)-g(a)\bigr).
\]
Thus evaluation at $a+v$ extends to $R_0$.

To choose the image of $z$, first evaluate the coefficients of $p$
at $a$. The resulting polynomial $p(a;Z)\in k[Z]$ is monic of
degree $\deg p>0$. Since $k$ is algebraically closed, choose $c\in k$
with $p(a;c)=0$. Evaluating $up+wp'=1$ at $(a;c)$ gives $w(a;c)p'(a;c)=1$.
The scalar part of $p(a+v;c)$ is $p(a;c)=0$, so $p(a+v;c)\in V$. For every
$\eta\in V$, the binomial theorem gives
\[
 p(a+v;c+\eta)=p(a+v;c)+p'(a+v;c)\eta
                =p(a+v;c)+p'(a;c)\eta.
\]
Indeed, $\eta^2=0$ and $(p'(a+v;c)-p'(a;c))\eta=0$ since $V^2=0$. Taking
$\eta=-p'(a;c)^{-1}p(a+v;c)$ gives $p(a+v;c+\eta)=0$.

By the presentation $R_0[Z]/(p)\xrightarrow{\sim}R$ above, the rule
$(h(z))\sigma=h(a+v;c+\eta)$ defines a $k$-algebra homomorphism
$\sigma\colon R\to C$ with $(t_i)\sigma=a_i+v_i$ and $(z)\sigma=c+\eta$.
For each $b\in\mathcal B$, both $b$ and $b^{-1}$ lie in $R$, so
$((b)\sigma)((b^{-1})\sigma)=1$. Thus $(b)\sigma$ is a unit.
\end{proof}

We shall change coefficients both to a residue field and to a larger
field. Generic freeness allows us to use the same localization for
both changes of coefficients.

\noindent\textbf{In the next two propositions,} $R$ is a
commutative integral domain, $B$ is an $R$-algebra, and the $B$-modules
$X,Y,Z$ under consideration are finitely generated over $R$.

\begin{proposition}\label{lem:upper-base}
Suppose that $R$ is noetherian and that $B$ is finitely generated
as an $R$-module.
For a $B$-module homomorphism
$u\colon Y\to Z$, put $Q=\Coker\bigl(\Hom_B(X,u)\bigr)$.

There is a nonzero $d\in R$ such that $\Hom_B(X,Y)$, $\Hom_B(X,Z)$
and $Q$ become free of finite rank over $R[1/d]$ after localization.
For every commutative $R$-algebra $T$ in which the image of $d$ is a
unit, there are natural $T$-module isomorphisms
\begin{align*}
 \Hom_B(X,U)\otimes_R T
 &\xrightarrow{\sim}\Hom_{B\otimes_R T}(X\otimes_R T,U\otimes_R T)
 \qquad(U=Y,Z),\\
 Q\otimes_R T
 &\xrightarrow{\sim}\Coker\bigl(\Hom_{B\otimes_R T}(X\otimes_R T,u\otimes_R 1_T)\bigr).
\end{align*}
For a $B$-module $M$, the action on $M\otimes_R T$ is given by
$(b\otimes t)(x\otimes t')=bx\otimes tt'$.
\end{proposition}

\begin{proof}
Any finite set of $R$-module generators of $X$ also generates $X$
over $B$. Choose an epimorphism $\pi\colon B^w\to X$ from such a set.
Since $B^w$ is finitely generated over the noetherian ring $R$,
$\Ker(\pi)$ is finitely generated over $R$. A finite set of generators
of $\Ker(\pi)$ therefore gives a presentation $B^v\to B^w\to X\to0$.
For $U=Y,Z$, applying $\Hom_B(-,U)$ gives an exact sequence
\[
 0\longrightarrow\Hom_B(X,U)\longrightarrow U^w
 \xrightarrow{\delta_U}U^v.
\]
The kernels, images and cokernels of $\delta_Y$ and $\delta_Z$, as
well as $Q$, are finitely generated over $R$. By generic freeness
\cite[\S~24, p.~185]{Matsumura}, each becomes free after localization
at a nonzero element of $R$. Take $d$ to be the product of these
elements. All the indicated modules are then free of finite rank over
$R[1/d]$ after localization.

For $U=Y,Z$, consider the short exact sequences
\[
\begin{gathered}
 0\longrightarrow\Ker(\delta_U)\longrightarrow U^w
   \longrightarrow\Image(\delta_U)\longrightarrow0,\\
 0\longrightarrow\Image(\delta_U)\longrightarrow U^v
   \longrightarrow\Coker(\delta_U)\longrightarrow0.
\end{gathered}
\]
After localization at $d$, the last term of each sequence is free,
so both sequences split. The structure homomorphism $R\to T$ extends
to $R[1/d]\to T$. Thus both sequences remain exact after tensoring
with $T$ over $R$. The first gives the kernel of
$U^w\otimes_R T\to(\Image(\delta_U))\otimes_R T$. The second shows that
$(\Image(\delta_U))\otimes_R T\to U^v\otimes_R T$ is injective.
Consequently, if $j\colon\Ker(\delta_U)\hookrightarrow U^w$ is the
inclusion, we have an isomorphism
\[
 (\Ker(\delta_U))\otimes_R T\xrightarrow{\sim}\Ker(\delta_U\otimes_R 1_T),
 \qquad x\otimes t\longmapsto(x)j\otimes t.
\]

Tensoring the presentation of $X$ with $T$ gives a presentation
of $X\otimes_R T$ by free $(B\otimes_R T)$-modules of ranks $v$ and
$w$. Applying $\Hom_{B\otimes_R T}(-,U\otimes_R T)$ gives the
exact sequence
\[
 0\longrightarrow\Hom_{B\otimes_R T}(X\otimes_R T,U\otimes_R T)
 \longrightarrow(U\otimes_R T)^w
 \xrightarrow{\delta_U\otimes_R 1_T}(U\otimes_R T)^v.
\]
Comparison of the two kernels gives the first isomorphism in the
statement, $f\otimes t\mapsto t(f\otimes_R1_T)$. Here
$t(f\otimes_R1_T)$ sends $x\otimes t'$ to $(x)f\otimes tt'$. For $f\colon X\to Y$ we have
\[
 (fu)\otimes_R1_T=(f\otimes_R1_T)(u\otimes_R1_T).
\]
By right exactness of tensor product, the isomorphism on cokernels
sends $[h]\otimes t$ to $[t(h\otimes_R1_T)]$ for $h\colon X\to Z$.
The brackets denote classes in the respective cokernels.
\end{proof}

In the proof of Proposition~\ref{prop:upper}, we must exclude
certain direct summands after reduction to the residue field.
The next proposition lifts a direct summand from the residue field
to the larger field.

\begin{proposition}\label{lem:upper-summand}
Let $L$ be a field containing $R$, and let $\mathfrak m$ be a maximal
ideal of $R$.
For $(M,N)=(X,Y)$ and $(Y,X)$, suppose the reduction homomorphism
\[
 \Hom_B(M,N)\longrightarrow
 \Hom_{B/\mathfrak m B}(M/\mathfrak m M,N/\mathfrak m N),\qquad
 f\longmapsto\bigl(x+\mathfrak m M\longmapsto(x)f+\mathfrak m N\bigr)
\]
is surjective. If $X/\mathfrak m X$ is isomorphic to a direct
summand of $Y/\mathfrak m Y$ over $B/\mathfrak m B$, then
$X\otimes_R L$ is isomorphic to a direct summand of $Y\otimes_R L$
over $B\otimes_R L$.
\end{proposition}

\begin{proof}
Choose a split monomorphism $\bar s\colon X/\mathfrak m X\to
Y/\mathfrak m Y$ and a retraction $\bar p$ over $B/\mathfrak m B$.
By hypothesis there are $B$-module homomorphisms $s\colon X\to Y$
and $p\colon Y\to X$ whose reductions are $\bar s$ and $\bar p$,
respectively. Since $\bar s\bar p=1_{X/\mathfrak m X}$, we have
$(x)(sp)-x\in\mathfrak m X$ for all $x\in X$. Hence
$X=\Image(sp)+\mathfrak m X$, and $H=\Coker(sp)$
satisfies $H=\mathfrak m H$.

Since $H$ is a quotient of $X$, the $R$-module $H$ is finitely
generated. Localizing at $\mathfrak m$ and applying Nakayama's lemma gives
$H\otimes_R R_{\mathfrak m}=0$. Every element of $R\setminus\mathfrak m$
is a unit in $L$, so $R\hookrightarrow L$ extends to $R_{\mathfrak m}\to L$.
Hence
\[
 H\otimes_R L
 \cong(H\otimes_R R_{\mathfrak m})\otimes_{R_{\mathfrak m}}L=0.
\]
Tensoring $X\xrightarrow{sp}X\to H\to0$ with $L$ shows that
$(sp)\otimes_R1_L$ is surjective. Since $X\otimes_R L$ has finite
dimension over $L$, $(sp)\otimes_R1_L$ is invertible. Thus
\[
 (p\otimes_R1_L)\bigl((sp)\otimes_R1_L\bigr)^{-1}
 \colon Y\otimes_R L\longrightarrow X\otimes_R L
\]
is a retraction of $s\otimes_R1_L$ over $B\otimes_R L$.
\end{proof}

To apply Lemma~\ref{lem:rank} in an induction on dimension, we need
endomorphisms that factor through a smaller module. When the
endomorphism ring is local, the images of the radical endomorphisms
generate a proper submodule.

\begin{proposition}\label{lem:upper-radical}
Let $B$ be an algebra over a field $F$, and let $M$ be a nonzero
$B$-module of finite dimension over $F$. Suppose $E=\End_B(M)$ is
local. Put $J=J(E)$ and
\[
 W=MJ=\sum_{h\in J}\Image(h).
\]
Let $i\colon W\hookrightarrow M$ be the inclusion. Then $W\subsetneq M$,
and the quotient map $E\to E/J$ gives a short exact sequence
\[
 0\longrightarrow\Hom_B(M,W)\xrightarrow{\Hom_B(M,i)}
 E\longrightarrow E/J\longrightarrow0.
\]
Moreover,
\begin{equation}\label{eq:upper-division}
 \dim_F W+\dim_F(E/J)\leq\dim_F M.
\end{equation}
\end{proposition}

\begin{proof}
As an $F$-subalgebra of $\End_F(M)$, the algebra $E$ is
finite-dimensional. Thus $J^a=0$ for some $a\geq1$. If
$M=MJ$, then $M=MJ^a=0$. Since $M\ne0$, the $B$-submodule $W$ is
proper.

Let $\pi\colon E\to E/J$ be the quotient homomorphism. We compare
the sequence in the statement with the exact bottom row in the
following diagram of $F$-modules.
\[
\begin{tikzcd}[column sep=2.2em,row sep=2.3em]
0 \arrow{r}
 & \Hom_B(M,W) \arrow{r}{\Hom_B(M,i)} \arrow{d}[swap]{\Phi}
 &[2.8em] E \arrow{r}{\pi} \arrow[equals]{d}
 & E/J \arrow{r} \arrow[equals]{d} & 0 \\
0 \arrow{r}
 & J \arrow[hook]{r}
 & E \arrow{r}[swap]{\pi}
 & E/J \arrow{r} & 0.
\end{tikzcd}
\]
For $g\colon M\to W$, set $(g)\Phi=gi$. Since $\Image(gi)\subseteq
W\subsetneq M$, the endomorphism $gi$ is not surjective and hence is
not a unit of the local algebra $E$. Therefore $gi\in J$.
Conversely, for $h\in J$ we have $\Image(h)\subseteq W$ by the definition
of $W$. Thus $h$ factors uniquely as $\widetilde h\,i$, where
$\widetilde h\colon M\to W$ sends $x$ to $(x)h$. The assignment $h\mapsto\widetilde h$ gives the
inverse of $\Phi$. The diagram commutes, so the top row is exact.

With endomorphisms acting on the right, $M$ is a right $E$-module.
Since $WE=M(JE)\subseteq W$ and $W=MJ$, the quotient $M/W$ is a
nonzero right $E/J$-module. Choose $x\in M\setminus W$.
The homomorphism of right $E/J$-modules
\[
 E/J\longrightarrow M/W,\qquad h+J\longmapsto(x)h+W,
\]
sends $1_M+J$ to $x+W\ne0$. The kernel is a proper right ideal of
the division algebra $E/J$ and hence is zero. Therefore
\[
 \dim_F(E/J)\leq\dim_F(M/W)=\dim_F M-\dim_F W.
\]
This proves \eqref{eq:upper-division}.
\end{proof}

We now apply the three propositions to the fields $k\subseteq L$
and the algebra $C$ of Lemma~\ref{lem:ring}. Before passing to $C$,
we write the action of $A$ over a subring of $L$.

Choose a $k$-basis $a_1,\ldots,a_s$ of $A$. Let $M$ be an $(L\otimes_k A)$-module with $L$-basis
$m_1,\ldots,m_n$. For $a\in A$, $\lambda\in L$ and $u\in M$,
write $au=(1_L\otimes a)u$ and $u\lambda=(\lambda\otimes1_A)u$.
Write
\[
 a_im_j=\sum_{h=1}^n m_hc_{hij}\qquad(c_{hij}\in L).
\]
For any $k$-subalgebra $R\subseteq L$ containing these coefficients,
put
\[
 M_{(R)}=\bigoplus_{j=1}^n m_jR\subseteq M.
\]
For each $i$, the element $a_i(\sum_j m_jr_j)=\sum_h m_h(\sum_j c_{hij}r_j)$ belongs to
$M_{(R)}$ for $r_1,\ldots,r_n\in R$. Thus $M_{(R)}$ is an
$A$-$R$-bimodule, free with basis $m_1,\ldots,m_n$ as a right
$R$-module. Multiplication gives an $(L\otimes_k A)$-module
isomorphism
\[
 M_{(R)}\otimes_R L\xrightarrow{\sim}M,\qquad
 u\otimes\lambda\longmapsto u\lambda.
\]

For a $k$-algebra homomorphism $\sigma\colon R\to C$, put
\[
 M^C=M_{(R)}\otimes_{R,\sigma}C,\qquad
 M_0=M^C\otimes_C\kappa,\qquad \kappa=C/V.
\]
The elements $m_j\otimes1_C$ form a right $C$-basis of $M^C$.
Using the scalar projection $\rho\colon C\to k$, we obtain
$A$-module isomorphisms
\[
 M^C/M^CV\xrightarrow{\sim}M_0
 \xrightarrow{\sim}M_{(R)}\otimes_{R,\sigma\rho}k.
\]
The quotient isomorphism sends $y+M^CV$ to $y\otimes\overline1$. Indeed,
tensoring $V\to C\to\kappa\to0$ with $M^C$ shows that the
surjection $y\mapsto y\otimes\overline1$ has kernel $M^CV$.
Tensor associativity and the isomorphism $\kappa\to k$,
$c+V\mapsto(c)\rho$, give the second map, which sends
$(u\otimes c)\otimes\overline1$ to $u\otimes(c)\rho$.

For a homomorphism $f\colon X\to Y$ of $(L\otimes_k A)$-modules,
choose $L$-bases of $X$ and $Y$. Suppose $R$ contains the
coefficients of the $A$-actions on $X$ and $Y$ and the entries of
the matrix of $f$ in these bases. Then $(X_{(R)})f\subseteq Y_{(R)}$.
Write $f_{(R)}\colon X_{(R)}\to Y_{(R)}$ for the resulting
$A$-$R$-bimodule homomorphism and put
\[
 f^C=f_{(R)}\otimes_R1_C\colon X^C\to Y^C,\qquad
 f_0=f^C\otimes_C1_\kappa\colon X_0\to Y_0.
\]
In the tensor products over $R$ with $k$, the formula for $f_0$ is
$x\otimes a\mapsto (x)f_{(R)}\otimes a$. Restriction and tensoring preserve identities and compositions.

\begin{proposition}\label{prop:upper}
Fix an integer $d_0\geq1$. Suppose there are only finitely many isomorphism classes
of indecomposable $A$-modules $N$ with $\dim_k N\leq d_0$ and
$N\not\cong\tau(N)$. Let $M$ be an $(L\otimes_k A)$-module with
$n=\dim_L M\leq d_0$, and fix an $L$-basis of $M$.

There are finite sets $\mathcal A\subseteq L$ and
$\mathcal B\subseteq L^\times$, with $\mathcal A$ containing the
coefficients $c_{hij}$ of the $A$-action in the chosen basis, such
that
\[
 \Delta(M^C)\leq n,\qquad M^C=M_{(R)}\otimes_{R,\sigma}C,
\]
whenever $R\subseteq L$ is a finitely generated $k$-subalgebra
containing $t_1,\ldots,t_m$ and $\mathcal A\cup\mathcal B$, and
$\sigma\colon R\to C$ is a $k$-algebra homomorphism satisfying
\[
 (t_i)\sigma=\lambda_i+v_i\quad(1\leq i\leq m),\qquad
 (b)\sigma\in C^\times\quad(b\in\mathcal B)
\]
for some $\lambda_1,\ldots,\lambda_m\in k$.

Moreover, for any finite sets
$\mathcal A\subseteq\mathcal A'\subseteq L$ and
$\mathcal B\subseteq\mathcal B'\subseteq L^\times$, one can choose
such $R$ and $\sigma$ with $\mathcal A'\cup\mathcal B'\subseteq R$
and $(b)\sigma\in C^\times$ for every $b\in\mathcal B'$.
\end{proposition}

\begin{proof}
We prove by induction on $n$ that there is a finitely generated
$k$-subalgebra $R_0\subseteq L$ containing $t_1,\ldots,t_m$ and all
$c_{hij}$ such that
\[
 \Delta(M_{(R_0)}\otimes_{R_0,\sigma} C)\leq n
 \qquad\text{for every $k$-algebra homomorphism $\sigma\colon R_0\to C$.}
\]
The same bound holds over every $k$-subalgebra $R\subseteq L$
containing $R_0$. Indeed, multiplication $M_{(R_0)}\otimes_{R_0} R\to M_{(R)}$
sends the $R$-basis $m_j\otimes1_R$ to the $R$-basis $m_j$, and
hence is an $A$-$R$-bimodule isomorphism. For every
$\sigma\colon R\to C$, tensor associativity gives
\[
 M_{(R)}\otimes_{R,\sigma}C\cong
 M_{(R_0)}\otimes_{R_0,\sigma|_{R_0}} C
\]
as $A$-$C$-bimodules. Isomorphic bimodules have the same value of
$\Delta$.

Start with the $k$-subalgebra $R_0$ of $L$ generated by
$t_1,\ldots,t_m$ and the coefficients $c_{hij}$. For $n=0$ we have $M^C=0$ and $\Delta(M^C)=0$. Suppose $n>0$ and first consider a
decomposition $f\colon M\xrightarrow{\sim}X\oplus Y$ of
$(L\otimes_k A)$-modules with $X,Y\ne0$. Choose
$L$-bases of $X$ and $Y$ and apply the induction hypothesis to these
two modules. Adjoin to $R_0$ generators of the resulting two
$k$-subalgebras and the matrix entries of $f$ and $f^{-1}$ in the
chosen bases. The restrictions of $f$ and $f^{-1}$ give mutually
inverse $A$-$R_0$-bimodule homomorphisms. Hence, for every
$\sigma\colon R_0\to C$, we have
\[
 M_{(R_0)}\xrightarrow[\sim]{f_{(R_0)}} X_{(R_0)}\oplus Y_{(R_0)},\qquad
 M^C\xrightarrow[\sim]{f^C} X^C\oplus Y^C.
\]
Subadditivity of $\Delta$ from Subsection~\ref{sec:rank} and the
induction hypothesis give
$\Delta(M^C)\leq\Delta(X^C)+\Delta(Y^C)\leq\dim_L X+\dim_L Y=n$.

Now suppose that $M$ is indecomposable as an $(L\otimes_k A)$-module.
Let $\mathcal N$ consist of one representative of each isomorphism class
of indecomposable $A$-modules $N$ such that $\dim_k N\leq n$ and
$N\not\cong\tau(N)$. The hypothesis on $d_0$ gives finiteness of $\mathcal N$.

Suppose that $f\colon M\to N\otimes_k L$ is an isomorphism of
$(L\otimes_k A)$-modules for some $N\in\mathcal N$. Choose a
$k$-basis $u_1,\ldots,u_n$ of $N$ and write
\[
 (m_j)f=\sum_i u_i\otimes\alpha_{ij},\qquad
 (u_i\otimes1_L)f^{-1}=\sum_j m_j\beta_{ji}.
\]
Adjoin all $\alpha_{ij}$ and $\beta_{ji}$ to $R_0$. The map
$N\otimes_k R_0\to N\otimes_k L$ given by the inclusion $R_0\subseteq L$
is injective with image the $R_0$-span of the $u_i\otimes1_L$.
The formulas for $f$ and $f^{-1}$ show that $f$ maps $M_{(R_0)}$ onto this span and that
$f^{-1}$ restricts to the inverse. For every $\sigma\colon R_0\to C$
we therefore have isomorphisms
\[
 M_{(R_0)}\xrightarrow[\sim]{f_{(R_0)}}N\otimes_k R_0,\qquad
 M^C\xrightarrow[\sim]{f_{(R_0)}\otimes_{R_0}1_C}
 (N\otimes_k R_0)\otimes_{R_0,\sigma} C
 \xrightarrow{\sim}N\otimes_k C.
\]
The first is an $A$-$R_0$-bimodule isomorphism. The last two are
$A$-$C$-bimodule isomorphisms, with the last map given by
$(z\otimes s)\otimes c\mapsto z\otimes(s)\sigma\,c$ and inverse
$z\otimes c\mapsto(z\otimes1_{R_0})\otimes c$.

For $\ell\in V^*$, tensor associativity now gives an isomorphism
$M^C\otimes_C\xi_\ell\cong N\otimes_k\xi_\ell$ of short exact
sequences of $A$-modules. Since $(n_2)p_\ell=\overline1$, the map
$N\otimes_k\kappa\to N\otimes_k N_\ell$ sending
$z\otimes\overline1$ to $z\otimes n_2$ is an $A$-module section of
$1_N\otimes_k p_\ell$. Thus $(\ell)d_{M^C}=0$ for all $\ell$,
and $\Delta(M^C)=0$.

Suppose finally that $M\not\cong N\otimes_k L$ as
$(L\otimes_k A)$-modules for every $N\in\mathcal N$. The algebra
$E=\End_{L\otimes_k A}(M)$ is local since $M$ is indecomposable.
Let $J$ be the Jacobson radical of $E$ and put
\[
 W=MJ=\sum_{h\in J}\Image(h),\qquad i\colon W\hookrightarrow M.
\]
By Proposition~\ref{lem:upper-radical}, $W$ is an
$(L\otimes_k A)$-submodule of dimension $r<n$. Choose an $L$-basis
$w_1,\ldots,w_r$ of $W$. The induction hypothesis gives a finitely
generated $k$-subalgebra $R_1\subseteq L$ for $W$ in this basis.
Write $(w_j)i=\sum_h m_h\gamma_{hj}$. Adjoin to $R_0$ a finite set of
$k$-algebra generators of $R_1$ and the entries $\gamma_{hj}$. Then
$W_{(R_0)}=\bigoplus_j w_jR_0\subseteq M_{(R_0)}$, so $i$ restricts to
$i_{(R_0)}\colon W_{(R_0)}\hookrightarrow M_{(R_0)}$. Multiplication
$W_{(R_1)}\otimes_{R_1}R_0\to W_{(R_0)}$ sends the basis $w_j\otimes1_{R_0}$
to $w_j$. Hence the induction hypothesis applied to $\sigma|_{R_1}$
gives $\Delta(W^C)\leq r$ for every $\sigma\colon R_0\to C$,
also after any further enlargement of $R_0$.

We apply Proposition~\ref{lem:upper-base} to the $R_0$-algebra
$R_0\otimes_k A$. The ring $R_0$ is a noetherian domain
since $R_0$ is a finitely generated $k$-subalgebra of a field.
The algebra $R_0\otimes_k A$ and the modules $M_{(R_0)}, W_{(R_0)}$
and $N\otimes_k R_0$ for $N\in\mathcal N$ are free of finite rank over
$R_0$. First take $X=Z=M_{(R_0)}$, $Y=W_{(R_0)}$ and $u=i_{(R_0)}$. For
each $N\in\mathcal N$ take also
\[
 (X,Y)=(N\otimes_k R_0,M_{(R_0)})\quad\text{and}\quad
 (X,Y)=(M_{(R_0)},N\otimes_k R_0),
\]
with $Z=0$ and $u\colon Y\to0$. Let $d$ be the product of the
nonzero elements supplied by all these applications. Localize at $d$
and extend the modules to $R_0[1/d]$ by taking the spans of the same
bases. The ring $R_0[1/d]\subseteq L$ is again a finitely generated
$k$-algebra. Write $R_0$ for this ring. Each factor of $d$ is now a
unit. For every commutative $R_0$-algebra $T$, the proposition gives
the natural $T$-module isomorphisms
\[
 \Hom_{R_0\otimes_k A}(X,U)\otimes_{R_0} T\xrightarrow{\sim}
 \Hom_{T\otimes_k A}(X\otimes_{R_0} T,U\otimes_{R_0} T)
 \qquad(U=Y,Z)
\]
for the modules $X,Y,Z$ listed above. Indeed, apply the proposition
over the original ring twice, first with the localization as target
and then with $T$ as target. Tensor associativity gives the displayed
isomorphism with $f\otimes t\mapsto t(f\otimes_{R_0}1_T)$.

In particular, put
\[
 Q=\Coker\!\left(
 \Hom_{R_0\otimes_k A}(M_{(R_0)},W_{(R_0)})
 \xrightarrow{\Hom_{R_0\otimes_k A}(M_{(R_0)},i_{(R_0)})}
 \End_{R_0\otimes_k A}(M_{(R_0)})\right).
\]
By the freeness and cokernel assertions of
Proposition~\ref{lem:upper-base} for $u=i_{(R_0)}$, the $R_0$-module
$Q$ is free of finite rank and, for every commutative $R_0$-algebra
$T$, there is an exact sequence
\[
\begin{aligned}
 \Hom_{T\otimes_k A}(M_{(R_0)}\otimes_{R_0} T,W_{(R_0)}\otimes_{R_0} T)
 &\longrightarrow \End_{T\otimes_k A}(M_{(R_0)}\otimes_{R_0} T)\\
 &\longrightarrow Q\otimes_{R_0} T\longrightarrow0,
\end{aligned}
\]
whose first map sends $g$ to $g(i_{(R_0)}\otimes_{R_0}1_T)$.

Fix a $k$-algebra homomorphism $\sigma\colon R_0\to C$. The composite
$\sigma\rho\colon R_0\to k$ makes $k$ an $R_0$-algebra. Since
$\sigma\rho$ restricts to the identity on $k$, the kernel
$\mathfrak m=\Ker(\sigma\rho)$ is maximal. The isomorphism $R_0/\mathfrak m\to k$
sends $s+\mathfrak m$ to $(s)(\sigma\rho)$. For $U=M,W$, write
\[
 \eta_U\colon U_{(R_0)}\otimes_{R_0} k\xrightarrow{\sim}U_0,\qquad
 u\otimes a\longmapsto (u\otimes1_C)\otimes(a+V).
\]
For each $N\in\mathcal N$ write
$\eta_N\colon(N\otimes_k R_0)\otimes_{R_0} k\to N$ for the isomorphism
$(z\otimes s)\otimes a\mapsto ((s)(\sigma\rho)a)z$.

For every $R_0$-module $P$, the map $P\to P\otimes_{R_0} k$ sending $x$ to
$x\otimes1$ is surjective since $P\otimes_{R_0} k\cong P/\mathfrak m P$.
Take $P$ to be each of the three Hom modules below and compose
with the base-change isomorphism for $T=k$. We use
$(R_0\otimes_k A)\otimes_{R_0} k\xrightarrow{\sim}A$,
$(s\otimes b)\otimes a\mapsto ((s)(\sigma\rho)a)b$. Composition with
$\eta_M$, $\eta_W$ and $\eta_N$ gives the epimorphisms
\begin{align}
 \Hom_{R_0\otimes_k A}(N\otimes_k R_0,M_{(R_0)})
 &\twoheadrightarrow\Hom_A(N,M_0),
 \label{eq:upper-reduction-N-to-M}\\
 \Hom_{R_0\otimes_k A}(M_{(R_0)},N\otimes_k R_0)
 &\twoheadrightarrow\Hom_A(M_0,N),
 \label{eq:upper-reduction-M-to-N}\\
 \Hom_{R_0\otimes_k A}(M_{(R_0)},W_{(R_0)})
 &\twoheadrightarrow\Hom_A(M_0,W_0).
 \label{eq:upper-reduction-M-to-W}
\end{align}
Writing $f_k=f\otimes_{R_0}1_k$, the maps send $f$ to
$\eta_N^{-1}f_k\eta_M$, $\eta_M^{-1}f_k\eta_N$ and
$\eta_M^{-1}f_k\eta_W$, respectively. The first two formulas hold for every $N\in\mathcal N$.

If $N\in\mathcal N$ were isomorphic to a direct summand of $M_0$
as an $A$-module, then \eqref{eq:upper-reduction-N-to-M} and
\eqref{eq:upper-reduction-M-to-N} would give the two surjections
required by Proposition~\ref{lem:upper-summand} at $\mathfrak m$. Hence
$N\otimes_k L$ would be isomorphic to a direct summand of $M$ as an
$(L\otimes_k A)$-module. Since $N\ne0$ and $M$ is indecomposable,
$M\cong N\otimes_k L$, contrary to the assumption $M\not\cong N\otimes_k L$ for all $N\in\mathcal N$.
The chosen basis of $M_{(R_0)}$ gives $\dim_k M_0=n$. Each indecomposable direct summand
$U$ of $M_0$ therefore has $\dim_k U\leq n$. If
$U\not\cong\tau(U)$, then $U$ is isomorphic to a member of
$\mathcal N$, which we have just excluded. Hence $U\cong\tau(U)$
as $A$-modules for every such summand. By additivity of $\tau$,
$M_0\cong\tau(M_0)$ as $A$-modules.

Put $i^C=i_{(R_0)}\otimes_{R_0}1_C$, $i_0=i^C\otimes_C1_\kappa$ and
\[
 J'=\{gi_0\mid g\in\Hom_A(M_0,W_0)\}\subseteq\End_A(M_0).
\]
As $(i_{(R_0)}\otimes_{R_0}1_k)\eta_M=\eta_W i_0$, taking $T=k$ in
the exact sequence for $Q$ gives a $k$-module isomorphism
$Q\otimes_{R_0} k\cong\End_A(M_0)/J'$. For $T=L$, the multiplication
isomorphisms $\mu_U\colon U_{(R_0)}\otimes_{R_0} L\xrightarrow{\sim}U$,
$u\otimes\lambda\mapsto u\lambda$, give the commutative diagram of $(L\otimes_k A)$-modules
\[
\begin{tikzcd}[column sep=large,row sep=large]
 W_{(R_0)}\otimes_{R_0} L
   \arrow{r}{i_{(R_0)}\otimes_{R_0}1_L} \arrow{d}[swap]{\mu_W}
 & M_{(R_0)}\otimes_{R_0} L \arrow{d}{\mu_M} \\
 W \arrow[hook]{r}{i} & M.
\end{tikzcd}
\]
Each path sends $w\otimes\lambda$ to $w\lambda$ in $M$. Thus
$Q\otimes_{R_0} L\cong\Coker\bigl(\Hom_{L\otimes_k A}(M,i)\bigr)$.
The image of $\Hom_{L\otimes_k A}(M,i)$ is $J$ by
Proposition~\ref{lem:upper-radical}. Since $Q$ is free of finite
rank over $R_0$, we obtain
\[
 \dim_k\bigl(\End_A(M_0)/J'\bigr)
 =\dim_k(Q\otimes_{R_0} k)=\operatorname{rank}_{R_0} Q
 =\dim_L(Q\otimes_{R_0} L)=\dim_L(E/J).
\]

For $h\in J'$, choose $g\colon M_0\to W_0$ with $h=gi_0$.
By \eqref{eq:upper-reduction-M-to-W} there is an
$(R_0\otimes_k A)$-module homomorphism $\widetilde g\colon
M_{(R_0)}\to W_{(R_0)}$ with
$(\widetilde g\otimes_{R_0}1_k)\eta_W=\eta_M g$. The $A$-$C$-bimodule
homomorphisms
\[
 M^C\xrightarrow{\widetilde g\otimes_{R_0}1_C}W^C
     \xrightarrow{i^C}M^C
\]
therefore reduce to $g$ and $i_0$, respectively. In particular,
every member of a $k$-basis of $J'$ has the factorization required
in Lemma~\ref{lem:rank}. The bimodules $M^C$ and $W^C$ are free of
finite rank over $C$ on the right, and $M_0\cong\tau(M_0)$. The
lemma, the induction bound for $W$ and \eqref{eq:upper-division}
now give
\[
\begin{aligned}
 \Delta(M^C)
 &\leq \Delta(W^C)+\dim_k\bigl(\End_A(M_0)/J'\bigr)\\
 &\leq \dim_L W+\dim_L(E/J)\leq n.
\end{aligned}
\]
The induction is complete.

Take $\mathcal A$ to consist of a finite set of $k$-algebra
generators of the resulting $R_0$ together with all $c_{hij}$, and
take $\mathcal B=\varnothing$. Every $R$ in the statement then
contains $R_0$, so the bound holds for every such $R$ and $\sigma$.
Applying Lemma~\ref{lem:ring} to $\mathcal A'$ and $\mathcal B'$
gives the last assertion.
\end{proof}

\section{Wild embeddings and homogeneity}\label{sec:wild}

Zhang and Xu stated the converse in Theorem~3 of their 2014 preprint
\cite{ZX}. In case~(III-4) of the proof of Proposition~5.5.6, the termination
argument requires a strict increase in the ordered pair $(\rho,\zeta)$.
Zhang and Xu cite Theorem~5.4.3(iv) for this step, but we have been
unable to deduce the required inequality from that theorem. We give a different proof using
self-extensions.

Let $T$ and $r$ be as in Theorem~\ref{thm:main}. We construct a simple
$S$-module $U$ of dimension $r$ and $r^2+1$ independent classes in
$\Ext_S^1(U,U)$. Their images under $T\otimes_S-$ will contradict the
upper bound of Proposition~\ref{prop:upper} if all but finitely many
indecomposable modules in each dimension up to $r^2$ are homogeneous.

\begin{proof}[Proof of Theorem~\ref{thm:main}]
Fix a right $S$-basis $z_1,\ldots,z_r$ of $T$. Suppose that in each
dimension $1\leq d\leq r^2$ there are only finitely many isomorphism classes of
indecomposable $A$-modules $X$ with $X\not\cong\tau_A(X)$. Taking
the union over $1\leq d\leq r^2$ gives the hypothesis of
Proposition~\ref{prop:upper} with $d_0=r^2$.

Put $m=r^2+1$, let $P=k[t_1,\ldots,t_m]$, and let $L$ be a
separable closure of the fraction field of $P$. Let $V$ have $k$-basis
$v_1,\ldots,v_m$. As in Lemma~\ref{lem:ring}, set $C=k\oplus V$ with
$V^2=0$ and put $\kappa=C/V$. Choose an $r\times r$ matrix
$(b_{ij})$ over $P$ as follows. Set $b_{i,i+1}=1$ for $1\leq i<r$.
Fill the remaining $r^2-r+1$ entries with $t_{r+1},\ldots,t_m$,
using each variable once. On the standard right $P$-basis
$e_1,\ldots,e_r$ of $P^r$, set
\[
 xe_j=e_jt_j,\qquad ye_j=\sum_{i=1}^r e_i b_{ij}
 \qquad(1\leq j\leq r).
\]
Extend $x$ and $y$ to endomorphisms of the right $P$-module $P^r$.
Since $S$ is free on $x,y$, these endomorphisms make $P^r$ an
$S$-$P$-bimodule.

Use the same matrices on the standard basis of $L^r$, also denoted by
$e_1,\ldots,e_r$. Then
\[
 P^r\otimes_P L\xrightarrow{\sim}L^r,\qquad
 \sum_{j=1}^r e_j\otimes c_j\longmapsto(c_1,\ldots,c_r),
\]
is an $S$-$L$-bimodule isomorphism. Here $s(u\otimes c)=(su)\otimes c$
for $s\in S$. Form the $(L\otimes_k A)$-module
\[
 M=T\otimes_S L^r,\qquad (\lambda\otimes a)(z\otimes u)=az\otimes u\lambda.
\]
The tensors $z_p\otimes e_j$, $1\leq p,j\leq r$, form an $L$-basis of
$M$, so $\dim_L M=r^2$.

Choose $\mathcal A\subseteq L$ and $\mathcal B\subseteq L^\times$
as in Proposition~\ref{prop:upper} for $M$ with this basis. If $r>1$,
enlarge $\mathcal B$ by $b_{r,1}$ and the elements $t_i-t_j$ with
$1\leq i<j\leq r$. All these elements are nonzero. By
Lemma~\ref{lem:ring}, there is a finitely generated $k$-subalgebra
$R$ of $L$ containing $P$ and $\mathcal A\cup\mathcal B$. We may
choose $R$ together with a $k$-algebra homomorphism $\sigma\colon R\to C$
such that
\[
 (t_\alpha)\sigma=\lambda_\alpha+v_\alpha
 \quad(1\leq\alpha\leq m),\qquad
 (b)\sigma\in C^\times\quad(b\in\mathcal B)
\]
for some $\lambda_1,\ldots,\lambda_m\in k$. Proposition~\ref{prop:upper}
now gives
\[
 \Delta(M^C)\leq r^2,\qquad M^C=M_{(R)}\otimes_{R,\sigma}C.
\]

The matrices of $x$ and $y$ have entries in $P$ and therefore preserve
$R^r=\bigoplus_{j=1}^r e_jR\subseteq L^r$.
For the inclusion $\iota\colon R^r\hookrightarrow L^r$, the
$A$-$R$-bimodule homomorphism
\[
 1_T\otimes_S\iota\colon T\otimes_S R^r\longrightarrow T\otimes_S L^r,
 \qquad z\otimes u\longmapsto z\otimes (u)\iota,
\]
takes the right $R$-basis consisting of the tensors $z_p\otimes e_j$
to the chosen $L$-basis of $M$. Hence $1_T\otimes_S\iota$ is injective
with image $M_{(R)}=\bigoplus_{p,j}(z_p\otimes e_j)R$. Thus $1_T\otimes_S\iota$ gives an $A$-$R$-bimodule isomorphism onto $M_{(R)}$.

Write $e_1^C,\ldots,e_r^C$ for the standard right $C$-basis of $C^r$.
Give $C^r$ the $S$-action determined by
\[
 xe_j^C=e_j^C(\lambda_j+v_j),\qquad
 ye_j^C=\sum_{i=1}^r e_i^C\bigl((b_{ij})\sigma\bigr)
 \qquad(1\leq j\leq r),
\]
where $x$ and $y$ act by endomorphisms of the right $C$-module.
On $R^r\otimes_{R,\sigma}C$, we have $s(u\otimes c)=(su)\otimes c$.
The right $C$-module isomorphism
\[
 R^r\otimes_{R,\sigma}C\xrightarrow{\sim}C^r,\qquad
 e_j\otimes c\longmapsto e_j^C c,
\]
commutes with $x$ and $y$ by the formulas above. Tensor associativity
therefore gives the $A$-$C$-bimodule isomorphism
\[
 T\otimes_S C^r\xrightarrow{\sim}M_{(R)}\otimes_{R,\sigma}C=M^C,
 \qquad z\otimes e_j^C c\longmapsto (z\otimes e_j)\otimes c.
\]

Put $U=C^r\otimes_C\kappa$ and $\bar e_j=e_j^C\otimes(1_C+V)$.
The isomorphism $\kappa\to k$, $(a+v)+V\mapsto a$, sends $1_C+V$
to $1$. Hence $\bar e_1,\ldots,\bar e_r$ form a $k$-basis of $U$.
Write $(b_{ij})\sigma=\beta_{ij}+w_{ij}$ with $\beta_{ij}\in k$
and $w_{ij}\in V$. Since $V\kappa=0$, the actions on $U$ are
\[
 x\bar e_j=\lambda_j\bar e_j,\qquad
 y\bar e_j=\sum_{i=1}^r\beta_{ij}\bar e_i.
\]

We claim that $U$ is simple. If $r=1$, then $U$ has dimension one.
Suppose $r>1$. The elements $(t_i-t_j)\sigma$ for $i<j$ and
$(b_{r,1})\sigma$ are units of $C$, so their images in $\kappa$ are nonzero. Hence
$\lambda_i\ne\lambda_j$ for $i\ne j$ and $\beta_{r,1}\ne0$.
The projection of $U$ onto $k\bar e_i$ along the other basis vectors
is the action of
\[
 \prod_{\substack{1\leq j\leq r\\j\ne i}}
 \frac{x-\lambda_j1_S}{\lambda_i-\lambda_j}\ \in S.
\]
Thus each projection preserves every $S$-submodule $U'$ of $U$.
If $U'\ne0$, choose $\sum_j c_j\bar e_j\in U'$ with $c_j\ne0$
for some $j$. Projection gives $c_j\bar e_j\in U'$, so $\bar e_j\in U'$.
Whenever $\bar e_{i+1}\in U'$, the projection of $y\bar e_{i+1}$
onto $k\bar e_i$ gives $\beta_{i,i+1}\bar e_i=\bar e_i\in U'$.
Starting at $j$, we obtain $\bar e_1\in U'$. The projection of
$y\bar e_1$ onto $k\bar e_r$ gives $\beta_{r,1}\bar e_r\in U'$.
Since $\beta_{r,1}\ne0$, we have $\bar e_r\in U'$. Starting now
at $r$, the implication $\bar e_{i+1}\in U'\Rightarrow\bar e_i\in U'$
gives every basis vector of $U$. Thus $U'=U$.

We next show that the homomorphism
\[
 d_{C^r}\colon V^*\longrightarrow\Ext^1_S(U,U),\qquad
 \ell\longmapsto[C^r\otimes_C\xi_\ell],
\]
from Subsection~\ref{sec:rank} is injective. If $(\ell)d_{C^r}=0$,
the sequence
\[
 0\longrightarrow U\xrightarrow{\,1_{C^r}\otimes_C j_\ell\,}
 C^r\otimes_C N_\ell\xrightarrow{\,1_{C^r}\otimes_C p_\ell\,}
 U\longrightarrow0
\]
splits in $\modu{S}$. Using the basis $n_1,n_2$ of
$N_\ell$ from Subsection~\ref{sec:rank}, take coordinates via the
$k$-module isomorphism
\[
 C^r\otimes_C N_\ell\xrightarrow{\sim}U\oplus U,\qquad
 e_j^C\otimes n_1\longmapsto(\bar e_j,0),\qquad
 e_j^C\otimes n_2\longmapsto(0,\bar e_j).
\]
In these coordinates $1_{C^r}\otimes_C p_\ell$ is the projection onto
the second factor. An $S$-module section $s_\ell$ therefore satisfies
$(\bar e_j)s_\ell=(u_j,\bar e_j)$ for some $u_j\in U$.

Compare the first coordinates in
$x((\bar e_j)s_\ell)=(x\bar e_j)s_\ell$. On $N_\ell$ we have
$Vn_1=0$ and $v_jn_2=((v_j)\ell)n_1$. Since
$xe_j^C=e_j^C(\lambda_j+v_j)$, the first coordinate on the left is
$xu_j+((v_j)\ell)\bar e_j$. The first coordinate on the right is $\lambda_j u_j$. Hence
\[
 (x-\lambda_j1_S)u_j=-((v_j)\ell)\bar e_j.
\]
The two sides lie in the complementary subspaces
$\bigoplus_{i\ne j}k\bar e_i$ and $k\bar e_j$, so both are zero.
Hence $(v_j)\ell=0$ for $1\leq j\leq r$. The $\lambda_j$-eigenspace
of $x$ is $k\bar e_j$, so $u_j=h_j\bar e_j$ for some $h_j\in k$.

The equality $y((\bar e_j)s_\ell)=(y\bar e_j)s_\ell$ gives
\[
 y(h_j\bar e_j,\bar e_j)
 =\sum_{i=1}^r\beta_{ij}(h_i\bar e_i,\bar e_i).
\]
Since $w_{ij}n_2=((w_{ij})\ell)n_1$, the first coordinates give
\[
 \sum_{i=1}^r\bigl(h_j\beta_{ij}+(w_{ij})\ell\bigr)\bar e_i
 =\sum_{i=1}^r\beta_{ij}h_i\bar e_i.
\]
Comparing the coefficients of $\bar e_i$ yields
\begin{equation}\label{eq:thm-a-section-y}
 (w_{ij})\ell=(h_i-h_j)\beta_{ij}
 \qquad(1\leq i,j\leq r).
\end{equation}
For $j=i+1$, the equality $b_{i,i+1}=1$ gives
$w_{i,i+1}=0$ and $\beta_{i,i+1}=1$. Thus $h_i=h_{i+1}$ for
$1\leq i<r$. All $h_i$ are equal, so \eqref{eq:thm-a-section-y}
gives $(w_{ij})\ell=0$ for every $i,j$.
Each $t_\alpha$ with $r<\alpha\leq m$ occurs as one of the $b_{ij}$.
For that entry $w_{ij}=v_\alpha$, since
$(t_\alpha)\sigma=\lambda_\alpha+v_\alpha$.
We have now shown $(v_\alpha)\ell=0$ both for $\alpha\le r$
from the action of $x$ and for $\alpha>r$ from the action of $y$.
Hence $\ell=0$, as required.

Since $T$ is free as a right $S$-module, tensoring Yoneda extensions
with $T$ over $S$ gives a $k$-module homomorphism
\[
 \Phi\colon\Ext^1_S(U,U)\longrightarrow\Ext^1_A(T\otimes_S U,T\otimes_S U),
 \qquad [\zeta]\longmapsto[T\otimes_S\zeta].
\]
A nonsplit self-extension of the simple module $U$ has indecomposable
middle term. Indeed, a decomposition would express that term as a sum of
two simple modules, and the extension would split. The representation
embedding preserves indecomposability of the middle term. Thus
$T\otimes_S\zeta$ is nonsplit whenever $\zeta$ is nonsplit, so
$\Phi$ is injective.

The $A$-$C$-bimodule isomorphism $T\otimes_S C^r\to M^C$ gives, for
each $C$-module $H$, an $A$-module isomorphism
\[
 \theta_H\colon T\otimes_S(C^r\otimes_C H)\xrightarrow{\sim}M^C\otimes_C H,
 \qquad z\otimes(e_j^C\otimes h)\longmapsto
 \bigl((z\otimes e_j)\otimes 1_C\bigr)\otimes h.
\]
The maps $\theta_H$ are natural in $H$. In particular, $\theta_\kappa$ maps
$T\otimes_S U$ isomorphically onto $M_0=M^C\otimes_C\kappa$.
Naturality for $j_\ell$ and $p_\ell$ gives an isomorphism from
$T\otimes_S(C^r\otimes_C\xi_\ell)$ to $M^C\otimes_C\xi_\ell$ with
$\theta_\kappa$ at both ends. To compare their classes, let
\[
 \Psi\colon\Ext^1_A(T\otimes_S U,T\otimes_S U)
 \xrightarrow{\sim}\Ext^1_A(M_0,M_0)
\]
be pullback along $\theta_\kappa^{-1}$ followed by pushout along
$\theta_\kappa$. If a class is represented by
$0\to T\otimes_S U\xrightarrow{\iota}Z\xrightarrow{\pi}T\otimes_S U\to0$,
the image of the extension class under $\Psi$ is represented by
\[
 0\longrightarrow M_0\xrightarrow{\theta_\kappa^{-1}\iota}
 Z\xrightarrow{\pi\theta_\kappa}M_0\longrightarrow0.
\]
Thus $d_{M^C}=d_{C^r}\Phi\Psi$. All three factors are injective, so
\[
 \Delta(M^C)=\dim_k\Image(d_{M^C})=\dim_k V^*=r^2+1.
\]
The equality $\Delta(M^C)=r^2+1$ contradicts the upper bound
$\Delta(M^C)\leq r^2$. Hence infinitely many
pairwise nonisomorphic indecomposable $A$-modules $X$ satisfy
$X\not\cong\tau_A(X)$ in some dimension $1\leq d\leq r^2$.
This proves the dimension bound.

For the last assertion, suppose first that $A$ is tame. Then $A$ is $\tau$-tame by Crawley-Boevey's homogeneity
theorem \cite[Theorem~D]{CB}. Conversely, suppose $A$ is $\tau$-tame.
A wild representation embedding would yield infinitely many pairwise
nonisomorphic nonhomogeneous indecomposable modules of the same dimension
by the dimension bound just proved. Thus $A$ is not wild. Drozd's
dichotomy~\cite{Drozd} implies that $A$ is tame.
\end{proof}

\begin{proof}[Proof of Theorem~\ref{thm:main2}]
If $A$ is tame, the assertion follows from
\cite[Corollary~E, p.~482]{CB}. Conversely, suppose the condition on
homogeneous tubes holds. Every vertex of a stable tube of rank one is
fixed by the Auslander--Reiten translation. Hence every indecomposable
$A$-module $X$ with $X\not\cong\tau_A(X)$ lies outside the homogeneous
tubes. In each dimension there are only finitely many such modules up
to isomorphism. Thus $A$ is $\tau$-tame, and Theorem~\ref{thm:main}
implies that $A$ is tame.
\end{proof}

Tameness passes to quotient algebras and to the opposite algebra.
Applying Theorem~\ref{thm:main} gives the following result.

\begin{corollary}\label{cor:descent}
Let $I$ be a two-sided ideal of $A$. If $A$ is $\tau$-tame,
then $A/I$ and $A\opp$ are $\tau$-tame.
\end{corollary}

\section{Extensions over radical quotients}\label{sec:ext}

\noindent\textbf{Throughout this section,} fix an integer $s\geq2$ and put
$J=J(A)$ and $B=A/J^s$. Regard each $B$-module as an $A$-module via
$A\twoheadrightarrow B$. The comparison of extension functors and the
dimension formula \eqref{eq:top-difference} hold over any field $k$.

We use this comparison to bound the last radical layers of quotients
of finite direct sums of $J$.

\subsection{Comparison of extensions}

An $A$-extension of two $B$-modules is a $B$-extension precisely
when $J^s$ annihilates its middle term. To express this condition in
terms of extension classes, consider the additive bifunctors
\[
 \Ext_B^1,\ \Ext_A^1\colon
 (\modu{B})\opp\times \modu{B}
 \longrightarrow \modu{k},
\]
where the arguments of $\Ext_A^1$ are restricted to $B$-modules.
For $B$-modules $X$ and $Y$, the homomorphism
$\gamma_{X,Y}\colon\Ext_B^1(X,Y)\to\Ext_A^1(X,Y)$ sends the class of
a $B$-extension
\[
 0\longrightarrow Y\xrightarrow{j}E\xrightarrow{p}X\longrightarrow0
\]
to the class of the same sequence with $A$-action $ae=(a+J^s)e$ on $E$.
Restriction of scalars preserves pullbacks and pushouts, so the maps
$\gamma_{X,Y}$ are natural in both variables. For fixed $X$, define
$C_X\colon \modu{B}\to \modu{k}$
by the exact sequence of covariant functors
\begin{equation}\label{eq:cokernel}
 \Ext_B^1(X,-)\xrightarrow{\gamma_{X,-}}\Ext_A^1(X,-)
 \xrightarrow{c_{X,-}}C_X\longrightarrow0.
\end{equation}
Thus $c_{X,Y}$ is the quotient map with kernel $\Image(\gamma_{X,Y})$.
The cokernels form an additive bifunctor
\[
 C\colon (\modu{B})\opp\times \modu{B}
 \longrightarrow \modu{k},\qquad (X,Y)\longmapsto C_X(Y).
\]
For $a\colon X'\to X$ and $u\colon Y\to Y'$, naturality of $\gamma$
gives
\[
 \bigl(([\xi])c_{X,Y}\bigr)C(a,u)
   =\bigl([u_*(a^*\xi)]\bigr)c_{X',Y'}.
\]
Here one pulls back $\xi$ along $a$ and then pushes out along $u$.

For a primitive idempotent $e_i$ of $A$, put $P_i=Ae_i$ and
$S_i=P_i/JP_i$. For $X=S_i$, write $\gamma=\gamma_{S_i,-}$ and
$c=c_{S_i,-}$. Let $a_i$ and $\iota_i$ denote the inclusions
\[
 J^sP_i\xrightarrow{a_i}JP_i\xrightarrow{\iota_i}P_i.
\]
For $M\in \modu{A}$, we use $\Hom_A(M,-)$ as a functor
from $\modu{B}$ to $\modu{k}$.
The component of $\Hom_A(a_i,-)$ at $Y$ is restriction to $J^sP_i$,
\[
 \Hom_A(JP_i,Y)\longrightarrow\Hom_A(J^sP_i,Y),
 \qquad f\longmapsto a_i f.
\]
Define $F_{i,s}\colon \modu{B}\to \modu{k}$
as the image of $\Hom_A(a_i,-)$ in the category of covariant
functors,
\begin{equation}\label{eq:restriction}
 F_{i,s}=\Image\!\left(
 \Hom_A(JP_i,-)\xrightarrow{\Hom_A(a_i,-)}\Hom_A(J^sP_i,-)
 \right).
\end{equation}
For a $B$-module homomorphism $u\colon Y\to Z$, we have
\[
 c_Y C_{S_i}(u)=\Ext_A^1(S_i,u)c_Z,\qquad
 (g)F_{i,s}(u)=gu\quad(g\in F_{i,s}(Y)).
\]
If $g=a_i f$, then $(J^sP_i)g=J^{s-1}((JP_i)f)\subseteq J^{s-1}Y$.
The composite $J^{s+1}P_i\hookrightarrow J^sP_i\xrightarrow{g}Y$
has image $J^s((JP_i)f)\subseteq J^sY=0$ and therefore equals zero.
Hence $g$ factors through the homomorphism of $A/J$-modules
\begin{equation}\label{eq:radical-layer-map}
 \bar g\colon J^sP_i/J^{s+1}P_i\longrightarrow J^{s-1}Y,
 \qquad x+J^{s+1}P_i\longmapsto(x)g.
\end{equation}

An extension of $S_i$ by $Y$ can be obtained by a pushout along a
homomorphism $JP_i\to Y$. The next proposition shows that the
restriction of this homomorphism to $J^sP_i$ depends only on the
extension class.

\begin{proposition}\label{prop:comparison}
Let $p_i\colon P_i\twoheadrightarrow S_i$ be the quotient map. For
$Y\in \modu{B}$ and $f\in\Hom_A(JP_i,Y)$, let $\xi_f$
be the pushout of
$0\to JP_i\xrightarrow{\iota_i}P_i\xrightarrow{p_i}S_i\to0$ along $f$.
There is a unique natural transformation
$\rho\colon\Ext_A^1(S_i,-)\to F_{i,s}$ such that
\begin{equation}\label{eq:restriction-ext}
 ([\xi_f])\rho_Y=a_i f
 \qquad(f\in\Hom_A(JP_i,Y)).
\end{equation}
For an $A$-extension $\xi\colon 0\to Y\xrightarrow{j}E\xrightarrow{p}S_i\to0$,
put $g=([\xi])\rho_Y$. Then $\Image(gj)=J^sE$, and
$[\xi]\in\Image(\gamma_Y)$ if and only if $J^sE=0$.

The sequence of functors
\[
 0\longrightarrow\Ext_B^1(S_i,-)\xrightarrow{\gamma}
 \Ext_A^1(S_i,-)\xrightarrow{\rho}F_{i,s}\longrightarrow0
\]
is exact. Consequently, there is a unique natural isomorphism
$\alpha\colon C_{S_i}\xrightarrow{\sim}F_{i,s}$ with $c\alpha=\rho$.
\end{proposition}

\begin{proof}
Fix a $B$-module $Y$ and apply $\Hom_A(-,Y)$ to the short exact
sequence with middle term $P_i$. Since $P_i$ is projective, we obtain
\[
 \Hom_A(P_i,Y)\xrightarrow{\Hom_A(\iota_i,Y)}
 \Hom_A(JP_i,Y)\xrightarrow{\delta_Y}\Ext_A^1(S_i,Y)
 \longrightarrow0.
\]
The first map sends $h$ to $\iota_i h$ and the connecting homomorphism
sends $f$ to $[\xi_f]$.
For $h\colon P_i\to Y$, we have $(J^sP_i)h\subseteq J^sY=0$,
so $a_i\iota_i h=0$. Hence
\[
 \Hom_A(\iota_i,Y)\,\Hom_A(a_i,Y)=0.
\]
Since $\delta_Y$ is surjective, exactness gives a unique
$k$-module homomorphism $\rho_Y$ satisfying \eqref{eq:restriction-ext}.
The image of $\rho_Y$ consists of all restrictions $a_i f$ and hence
equals $F_{i,s}(Y)$.

Let $u\colon Y\to Z$ be a $B$-module homomorphism. Pushing out
$\xi_f$ along $u$ gives an extension representing $[\xi_{fu}]$. Thus
\[
 ([\xi_f])\Ext_A^1(S_i,u)\rho_Z
 =a_i(fu)=(a_i f)u=([\xi_f])\rho_Y F_{i,s}(u).
\]
Every extension class has the form $[\xi_f]$, so $\rho$ is natural.

For the given extension $\xi$, use projectivity of $P_i$ to choose
$t\colon P_i\to E$ with $tp=p_i$. Since $(\iota_i t)p=\iota_i p_i=0$
and $j$ is a kernel of $p$, there is a unique $f\colon JP_i\to Y$
with $fj=\iota_i t$. These maps give a commutative diagram in
$\modu{A}$,
\[
\begin{tikzcd}[column sep=large, row sep=2em]
0 \arrow{r}
 & JP_i \arrow[hook]{r}{\iota_i} \arrow{d}[swap]{f}
 & P_i \arrow{r}{p_i} \arrow{d}{t}
 & S_i \arrow{r} \arrow[equals]{d} & 0 \\
0 \arrow{r}
 & Y \arrow[hook]{r}{j}
 & E \arrow{r}[swap]{p}
 & S_i \arrow{r} & 0.
\end{tikzcd}
\]
By the pushout property, there is a morphism $\xi_f\to\xi$ whose
components on $Y$ and $S_i$ are the identity. Thus $[\xi]=[\xi_f]$,
and \eqref{eq:restriction-ext} gives $g=a_i f$ and $gj=a_i\iota_i t$.
Since $tp=p_i$ is surjective and $\Ker(p)=\Image(j)$, we have
$E=\Image(j)+\Image(t)$. Therefore
\[
 J^sE=(J^sY)j+(J^sP_i)t=(J^sP_i)t=\Image(gj).
\]

Any $A$-module section of a $B$-extension is also a $B$-module
homomorphism. Hence $\gamma_Y$ is injective. Since $j$ is injective,
$g=0$ if and only if $\Image(gj)=0$, that is, $J^sE=0$. In this case
the action on $E$ factors through $B$, so $\xi$ is a $B$-extension.
Conversely, if $[\xi]\in\Image(\gamma_Y)$, then $E$ is isomorphic over
$A$ to the middle term of a $B$-extension. Hence $J^sE=0$.
Thus $\Ker(\rho_Y)=\Image(\gamma_Y)$. Since $\gamma_Y$ is injective and
$\rho_Y$ is surjective, the sequence is exact at every $Y$.
Now \eqref{eq:cokernel} gives the unique natural isomorphism $\alpha$
with $c\alpha=\rho$.
\end{proof}

\subsection{Bounds on the cokernel}\label{sec:cokerbounds}

Put $\Gamma=A/J$. Every simple $A$-module is a direct summand of the
regular $\Gamma$-module. We use the $\Gamma$-action on $C_\Gamma(Y)$
to compute $C_{S_i}(Y)$ for each $i$.
We regard $\Gamma$ as a module over $A$ or $B$ via the respective
quotient map.
For $\Lambda\in\{A,B\}$ and a $\Lambda$-module $Z$, evaluation at
$1_\Gamma$ gives a natural isomorphism from $\Hom_\Lambda(\Gamma,Z)$
to the socle of $Z$.
Indeed, $J(B)=J/J^s$, so in both cases the socle consists of the
elements annihilated by $J$. Thus $\Ext^1_\Lambda(\Gamma,-)$
is the first right derived functor of the socle functor.
For $t\in\Gamma$, let $r_t\colon\Gamma\to\Gamma$ be right
multiplication by $t$. Pullback along $r_t$ defines a left
$\Gamma$-action on $C_\Gamma(Y)$ by
\[
 t w=(w)C(r_t,1_Y)\qquad(w\in C_\Gamma(Y)).
\]
Write $\Head(Z)=Z/J(\Lambda)Z$ and
$\tau_\Lambda^-(Z)=\Tr_{\Lambda\opp}(DZ)$ for a $\Lambda$-module
$Z$. For an $A$-module $X$, put
$\ann_X(J^s)=\{x\in X\mid ax=0\text{ for every }a\in J^s\}$.

The injective envelopes of $Y$ over $A$ and $B$ give two descriptions
of extensions of $\Gamma$ by $Y$. Comparing the two descriptions gives
a formula for $\dim_k C_{S_i}(Y)$.

\begin{lemma}\label{lem:fixed}
For $X\in\modu{\Gamma}$ and $Y\in \modu{B}$,
there is an isomorphism of $k$-modules
\[
 C_X(Y)\xrightarrow{\sim}\Hom_\Gamma(X,C_\Gamma(Y))
\]
natural in $X$ and $Y$.

Choose an injective envelope $j_A\colon Y\hookrightarrow I_A(Y)$ over
$A$, and put
\[
 I_B(Y)=\ann_{I_A(Y)}(J^s).
\]
Then $j_B\colon Y\to I_B(Y)$, $y\mapsto(y)j_A$, is an injective
envelope over $B$. Let $q_A\colon I_A(Y)\twoheadrightarrow
Q_A=I_A(Y)/\Image(j_A)$ be the quotient map and put $Q_B=(I_B(Y))q_A$.
For $\Lambda\in\{A,B\}$, let $U_\Lambda$ be the socle of $Q_\Lambda$.
Then $U_B=U_A\cap Q_B$, and
\[
 C_\Gamma(Y)\cong U_A/U_B,\qquad
 C_X(Y)\cong\Hom_\Gamma(X,U_A/U_B).
\]
The first isomorphism respects the $\Gamma$-action, and the second is an isomorphism of $k$-modules.

Let $d_i=\dim_k\End_A(S_i)$ and write $[H:S_i]$ for the
multiplicity of $S_i$ as a composition factor of $H$. Then
\begin{equation}\label{eq:top-difference}
 \dim_k C_{S_i}(Y)=d_i\bigl(
 [\Head(\tau_A^-(Y)):S_i]
 -[\Head(\tau_B^-(Y)):S_i]\bigr).
\end{equation}
Over an algebraically closed field, $d_i=1$. If $A$ is also
$\tau$-tame, then for each $d\geq1$ only finitely many isomorphism
classes of indecomposable $B$-modules $Y$ of dimension $d$ satisfy
$C_{S_i}(Y)\ne0$, or equivalently $F_{i,s}(Y)\ne0$.
\end{lemma}

\begin{proof}
For $x\in X$, let $\ell_x\colon\Gamma\to X$ send $a$ to $ax$.
Define
\[
 C_X(Y)\longrightarrow\Hom_\Gamma(X,C_\Gamma(Y)),\qquad
 w\longmapsto\bigl(x\longmapsto(w)C(\ell_x,1_Y)\bigr).
\]
Since $\ell_{tx}=r_t\ell_x$, the map $x\mapsto(w)C(\ell_x,1_Y)$
commutes with the action of $\Gamma$. Additivity follows from
$\ell_{x+x'}=\ell_x+\ell_{x'}$. For a $\Gamma$-module homomorphism
$h\colon X\to X'$, the identity $\ell_x h=\ell_{(x)h}$ gives
naturality in $X$. Naturality in $Y$ follows from bifunctoriality
of $C$. When $X=\Gamma$, evaluation at $1_\Gamma$ is the inverse.
Since $\Gamma$ is semisimple, every $X\in\modu{\Gamma}$
is a direct summand of some $\Gamma^{\oplus n}$. Additivity now gives
the isomorphism for each $X$.

Since $J^sY=0$, the image of $j_A$ lies in $\ann_{I_A(Y)}(J^s)$.
The quotient-ring formula of Sharpe and V\'amos, recalled in
\cite[Proposition~2.4]{Crivei}, gives the asserted injective envelope
$j_B$ over $B$. Write
$b\colon I_B(Y)\hookrightarrow I_A(Y)$ and
$v\colon Q_B\hookrightarrow Q_A$ for the inclusions. The map
$q_B\colon I_B(Y)\twoheadrightarrow Q_B$, $x\mapsto(x)q_A$,
satisfies $q_Bv=bq_A$, and $j_Bb=j_A$. Thus for
$\Lambda\in\{A,B\}$ we have a short exact sequence
\[
 0\longrightarrow Y\xrightarrow{j_\Lambda}I_\Lambda(Y)
 \xrightarrow{q_\Lambda}Q_\Lambda\longrightarrow0.
\]
In particular, $Q_\Lambda$ is the first cosyzygy of $Y$ over
$\Lambda$. Since $J(B)=J/J^s$, an element of $Q_B$ lies in $U_B$
exactly when the element lies in $U_A$. Hence $U_B=U_A\cap Q_B$.

For $z\in U_\Lambda$, the map $\ell_z\colon\Gamma\to Q_\Lambda$,
$a\mapsto az$, is well defined since $Jz=0$.
Applying $\Hom_\Lambda(\Gamma,-)$ to the short exact sequence
gives the exact row in the following diagram of $k$-modules.
\[
\begin{tikzcd}[column sep=1.6em, row sep=1.7em]
 & U_\Lambda \arrow{d}[swap]{\sim} \arrow{dr}{\delta_\Lambda} & & \\
 \Hom_\Lambda(\Gamma,I_\Lambda(Y)) \arrow{r}{0}
 & \Hom_\Lambda(\Gamma,Q_\Lambda) \arrow{r}[swap]{\partial_\Lambda}
 & \Ext^1_\Lambda(\Gamma,Y) \arrow{r} & 0.
\end{tikzcd}
\]
The vertical isomorphism sends $z$ to $\ell_z$. Every homomorphism
$f\colon\Gamma\to I_\Lambda(Y)$ has semisimple image. Since
$j_\Lambda$ is an essential embedding, this image lies in $\Image(j_\Lambda)$.
Thus the first horizontal arrow sends $f$ to $f q_\Lambda=0$.
The next term $\Ext^1_\Lambda(\Gamma,I_\Lambda(Y))$ is zero since
$I_\Lambda(Y)$ is injective.
Consequently the connecting homomorphism $\partial_\Lambda$ is an
isomorphism, and so is the composite
\[
 \delta_\Lambda\colon U_\Lambda\xrightarrow{\sim}\Ext^1_\Lambda(\Gamma,Y),
 \quad z\longmapsto
 [\ell_z^*(0\to Y\xrightarrow{j_\Lambda}I_\Lambda(Y)
                      \xrightarrow{q_\Lambda}Q_\Lambda\to0)].
\]
The composite commutes with the $\Gamma$-actions since
$\ell_{tz}=r_t\ell_z$ and $t$ acts on extensions by pullback
along $r_t$.

Fix $z\in U_B$. The pullbacks representing $(z)\delta_B$ and
$((z)v)\delta_A$ have middle terms
\[
\begin{aligned}
 E_B &= \{(a,x)\in\Gamma\oplus I_B(Y)\mid (x)q_B=az\},\\
 E_A &= \{(a,x)\in\Gamma\oplus I_A(Y)\mid (x)q_A=(az)v\}.
\end{aligned}
\]
For $(a,x)\in E_B$, we have $((x)b)q_A=((x)q_B)v=(az)v$.
Hence $h\colon E_B\to E_A$, $(a,x)\mapsto(a,(x)b)$, is an
$A$-module homomorphism. In each extension the injection
$\iota_\Lambda$ sends $y$ to $(0,(y)j_\Lambda)$ and the projection
$\pi_\Lambda$ sends $(a,x)$ to $a$. Thus $j_Bb=j_A$ gives
$\iota_B h=\iota_A$, and $h\pi_A=\pi_B$.
Thus $(1_Y,h,1_\Gamma)$ is a morphism of $A$-extensions.
Consequently $((z)\delta_B)\gamma_{\Gamma,Y}=((z)v)\delta_A$.
That is, the following diagram commutes in $\modu{k}$.
\[
\begin{tikzcd}[column sep=large, row sep=2em]
 U_B \arrow{r}{\delta_B} \arrow[hook]{d}[swap]{v|_{U_B}}
 & \Ext_B^1(\Gamma,Y) \arrow{d}{\gamma_{\Gamma,Y}} \\
 U_A \arrow{r}[swap]{\delta_A}
 & \Ext_A^1(\Gamma,Y).
\end{tikzcd}
\]
Taking cokernels gives the $\Gamma$-module isomorphism
\[
 C_\Gamma(Y)\xrightarrow{\sim}U_A/U_B,\qquad
 ((z)\delta_A)c_{\Gamma,Y}\longmapsto z+U_B.
\]
Applying $\Hom_\Gamma(X,-)$ and using the first isomorphism of the
lemma gives the description of $C_X(Y)$.

For the dimension formula, choose an injective envelope
$e\colon Q_\Lambda\hookrightarrow I^1_\Lambda(Y)$. Since $e$ is
essential, $(U_\Lambda)e$ is the socle of $I^1_\Lambda(Y)$.
Thus $I_\Lambda(S_i)$ occurs in $I^1_\Lambda(Y)$ with multiplicity
$[U_\Lambda:S_i]$. The multiplicity formula for a minimal injective resolution
\cite[Section~1.8, Lemma~1]{CBRes} now gives
\[
 \dim_k\Ext^1_\Lambda(S_i,Y)=d_i[U_\Lambda:S_i].
\]
Now $0\to Y\xrightarrow{j_\Lambda}I_\Lambda(Y)
\xrightarrow{q_\Lambda e} I^1_\Lambda(Y)$ is a minimal injective
copresentation. By \cite[Chapter~IV, Proposition~1.12]{ARS},
the top of $\tau_\Lambda^-(Y)$ is isomorphic to the socle of
$I^1_\Lambda(Y)$ when $Y$ is indecomposable and noninjective.
For an injective summand $Y'$ of $Y$, both $\tau_\Lambda^-(Y')$
and $I^1_\Lambda(Y')$ are zero. Additivity therefore gives this
description of the top for every $Y$. Hence
\[
 \dim_k\Ext^1_\Lambda(S_i,Y)=
 d_i[\Head(\tau_\Lambda^-(Y)):S_i].
\]
By Proposition~\ref{prop:comparison},
\[
 \dim_k C_{S_i}(Y)=\dim_k\Ext_A^1(S_i,Y)-\dim_k\Ext_B^1(S_i,Y).
\]
Substituting the formula for $\Lambda=A$ and $\Lambda=B$ gives
\eqref{eq:top-difference}.

For the last assertion, Schur's lemma gives $d_i=1$ over an
algebraically closed field. If $A$ is $\tau$-tame, so is $B$ by
Corollary~\ref{cor:descent}. Restriction of scalars along
$A\twoheadrightarrow B$ is fully faithful and preserves dimension and
indecomposability. Hence, in each dimension, all but finitely many
isomorphism classes of indecomposable $B$-modules are homogeneous
over both $A$ and $B$. For each such $Y$,
\[
 \tau_A^-(Y)\cong Y\cong\tau_B^-(Y)
\]
by \cite[Chapter~IV, Proposition~1.10]{ARS}. Since $J(B)Y=JY$, the
two multiplicities in \eqref{eq:top-difference} are equal.
Proposition~\ref{prop:comparison} gives the assertion for $F_{i,s}(Y)$.
\end{proof}

\noindent\textbf{For the rest of this section,} take $k$ algebraically
closed and let $Y\in\modu{B}$ be indecomposable. When $A$ is
$\tau$-tame, Lemma~\ref{lem:fixed} allows finitely many exceptional
modules in each dimension. To bound $\dim_k C_{S_i}(Y)$ for those
modules as well, we construct a representation embedding from five
independent restrictions.

\begin{lemma}\label{lem:four}
If there is a monomorphism of $k$-modules
$k^5\hookrightarrow F_{i,s}(Y)$, then there is an
$A$-$S$-bimodule $T$, free as a right $S$-module of rank
$2\dim_k Y+5\dim_k S_i$, such that $T\otimes_S-$ is a representation
embedding. Consequently, if $A$ is $\tau$-tame, there is an epimorphism
of $k$-modules $k^4\twoheadrightarrow F_{i,s}(Y)$, and
\[
 \dim_k C_{S_i}(Y)=\dim_k F_{i,s}(Y)\leq4.
\]
\end{lemma}

\begin{proof}
We first choose five elements of $F_{i,s}(Y)$. The choice of five
homomorphisms $JP_i\to Y$ with the prescribed restrictions will
give the bimodule $T$ by a pushout.

Put $F=F_{i,s}(Y)$ and $R=\End_B(Y)$. By \eqref{eq:restriction},
postcomposition makes $F$ a right $R$-module. Since $Y$ is
indecomposable, $R$ is local. The division algebra $R/J(R)$ is
finite-dimensional over the algebraically closed field $k$. Hence the
homomorphism
\[
 k\longrightarrow R/J(R),\qquad \lambda\longmapsto\lambda 1_Y+J(R)
\]
is an isomorphism, where $\lambda 1_Y\colon Y\to Y$ sends $y$ to
$\lambda y$. For each $r\in R$, let $\bar r\in k$ be the unique scalar
satisfying
\[
 r+J(R)=\bar r\,1_Y+J(R).
\]

Suppose that $\dim_k F\geq5$. Since $R/J(R)\cong k$, every simple
right $R$-module has dimension one over $k$. Thus a composition
series of the right $R$-module $F$ has $\dim_k F$ factors. Choose
$0=F_0\subset F_1\subset\cdots\subset F_5$ as the first five inclusions
in a composition series of $F$.
For $1\leq j\leq5$, choose $f_j\in F_j\setminus F_{j-1}$ and
$\widehat f_j\colon JP_i\to Y$ such that $a_i\widehat f_j=f_j$.
Since each $F_j/F_{j-1}$ has dimension one, $f_1,\ldots,f_j$ form
a $k$-basis of $F_j$. The action of $R$ on $F_j/F_{j-1}$ factors
through $R/J(R)$. For all $r\in R$ and $1\leq j\leq5$, we have
\[
 f_jr+F_{j-1}=f_j(\bar r\,1_Y)+F_{j-1}.
\]

For an $S$-module $V$ of arbitrary dimension, define five monomorphisms of $k$-modules
$a_j^V\colon V\to V^{\oplus2}$ by
\[
 (v)a_1^V=(v,0),\quad (v)a_2^V=(0,v),\quad
 (v)a_3^V=(v,v),\quad (v)a_4^V=(v,xv),\quad (v)a_5^V=(v,yv).
\]
These are the five maps in the representation embedding for the
five-subspace quiver described by Crawley-Boevey
\cite[Section~6.2, Example~F(b)]{CBNotes}.
For $V=S$, each $a_j^S$ is a homomorphism of right $S$-modules. Define an
$A$-$S$-bimodule homomorphism
\[
 \Phi_S\colon JP_i\otimes_k S^{\oplus5}\longrightarrow
 Y\otimes_k S^{\oplus2},\qquad
 z\otimes(s_1,\ldots,s_5)\longmapsto
 \sum_{j=1}^5 (z)\widehat f_j\otimes(s_j)a_j^S.
\]
Tensor the exact sequence
$0\to JP_i\xrightarrow{\iota_i}P_i\xrightarrow{p_i}S_i\to0$
over $k$ with $S^{\oplus5}$. Pushing out along $\Phi_S$ gives a
commutative diagram of $A$-$S$-bimodules with exact rows
\[
\begin{tikzcd}[column sep=2em,row sep=large]
 0 \arrow{r}
 & JP_i\otimes_k S^{\oplus5}
     \arrow[hook]{r}{\iota_i\otimes 1}
     \arrow{d}[swap]{\Phi_S}
 & P_i\otimes_k S^{\oplus5}
     \arrow{r}{p_i\otimes 1}
     \arrow{d}{\lambda_S}
 & S_i\otimes_k S^{\oplus5}
     \arrow{r} \arrow[equals]{d}
 & 0 \\
 0 \arrow{r}
 & Y\otimes_k S^{\oplus2} \arrow[hook]{r}[swap]{j_S}
 & T \arrow{r}[swap]{q_S}
 & S_i\otimes_k S^{\oplus5} \arrow{r}
 & 0.
\end{tikzcd}
\]
The bottom row splits over $S$ because $S_i\otimes_k S^{\oplus5}$
is free as a right $S$-module. Thus $T$ is free as a right $S$-module
of rank $2\dim_k Y+5\dim_k S_i$, and $E=T\otimes_S-$ is exact.

Let $V$ be a finite-dimensional $S$-module and put
$L_V=Y\otimes_k V^{\oplus2}$. For each $A$-module $X$ and $m\geq1$,
there is an $A$-module isomorphism
\[
 (X\otimes_k S^{\oplus m})\otimes_S V
 \xrightarrow{\sim} X\otimes_k V^{\oplus m},\qquad
 (z\otimes(s_1,\ldots,s_m))\otimes v
 \longmapsto z\otimes(s_1v,\ldots,s_mv).
\]
Tensor the diagram over $S$ with $V$ and apply these isomorphisms.
Write $\Phi_V,\lambda_V,j_V,q_V$ for the resulting maps and $\xi_V$
for the bottom row in
\begin{equation}\label{eq:four-extension}
\begin{tikzcd}[column sep=2em,row sep=large]
 0 \arrow{r}
 & JP_i\otimes_k V^{\oplus5}
     \arrow[hook]{r}{\iota_i\otimes 1}
     \arrow{d}[swap]{\Phi_V}
 & P_i\otimes_k V^{\oplus5}
     \arrow{r}{p_i\otimes 1}
     \arrow{d}{\lambda_V}
 & S_i\otimes_k V^{\oplus5}
     \arrow{r} \arrow[equals]{d}
 & 0 \\
 0 \arrow{r}
 & L_V \arrow[hook]{r}[swap]{j_V}
 & E(V) \arrow{r}[swap]{q_V}
 & S_i\otimes_k V^{\oplus5} \arrow{r}
 & 0.
\end{tikzcd}
\end{equation}
The rows are exact since the original rows split as sequences of
right $S$-modules. The square with vertical maps $\Phi_V$ and
$\lambda_V$ is a pushout since $-\otimes_S V$ preserves pushouts. For
$w=(v_1,\ldots,v_5)\in V^{\oplus5}$ and $z\in JP_i$, we have
\[
 (z\otimes w)\Phi_V=\sum_{j=1}^5(z)\widehat f_j\otimes(v_j)a_j^V.
\]

We next compute the submodule of $E(V)$ annihilated by $J^s$. For any
finite-dimensional $k$-module $U$, define
\[
 \beta_U\colon F\otimes_k U\longrightarrow F_{i,s}(Y\otimes_k U),
 \qquad f\otimes u\longmapsto\bigl(z\mapsto(z)f\otimes u\bigr).
\]
The displayed map on $J^sP_i$ extends to $JP_i$ since $f$ does.
Thus $\beta_U$ takes values in $F_{i,s}(Y\otimes_k U)$. For a
homomorphism of $k$-modules $l\colon U\to U'$, we have
$\beta_U F_{i,s}(1_Y\otimes l)=(1_F\otimes l)\beta_{U'}$.
For $U=k$, the map $\beta_k$ is an isomorphism. A basis of $U$
and additivity of $F_{i,s}$ now show that $\beta_U$ is an isomorphism.

For $w\in V^{\oplus5}$, let $p_w\colon S_i\to S_i\otimes_k V^{\oplus5}$
send $c$ to $c\otimes w$. The pullback of $\xi_V$ along $p_w$
has the form
\[
 p_w^*\xi_V\colon\quad
 0\longrightarrow L_V
 \xrightarrow{\ y\mapsto(0,(y)j_V)\ }E_w
 \xrightarrow{\ (c,e)\mapsto c\ }S_i\longrightarrow0,
\]
where $E_w=\{(c,e)\in S_i\oplus E(V)\mid(c)p_w=(e)q_V\}$.
Since $J^sL_V=0$, Proposition~\ref{prop:comparison} gives the component
$\rho_{L_V}\colon\Ext_A^1(S_i,L_V)\to F_{i,s}(L_V)$. Define
$\Theta_V\colon V^{\oplus5}\to F_{i,s}(L_V)$ by
\[
 (w)\Theta_V
 =([\xi_V])\Ext_A^1(p_w,L_V)\rho_{L_V}
 =([p_w^*\xi_V])\rho_{L_V}.
\]
For $w=(v_1,\ldots,v_5)$, put
\[
 t_w\colon P_i\to E(V),\quad z\mapsto(z\otimes w)\lambda_V,
 \qquad
 d_w\colon JP_i\to L_V,\quad z\mapsto(z\otimes w)\Phi_V.
\]
The right square in \eqref{eq:four-extension} gives
\[
 ((z)t_w)q_V=(z\otimes w)(p_i\otimes 1)=((z)p_i)p_w
 \qquad (z\in P_i).
\]
The left square gives
\[
 ((z)\iota_i)t_w=(z\otimes w)\Phi_Vj_V=((z)d_w)j_V
 \qquad (z\in JP_i).
\]
Hence $t_wq_V=p_ip_w$ and $\iota_i t_w=d_wj_V$.
Thus $z\mapsto((z)p_i,(z)t_w)$ defines a lift $P_i\to E_w$ of $p_i$.
On $JP_i$ the lift sends $z$ to $(0,((z)d_w)j_V)$.
Together with $1_{L_V}$, the lift gives a morphism from the pushout extension
$\xi_{d_w}$ to $p_w^*\xi_V$ with the identity on each end term.
Hence $[p_w^*\xi_V]=[\xi_{d_w}]$. By \eqref{eq:restriction-ext},
$(w)\Theta_V=a_i d_w$. In particular, for $z\in J^sP_i$,
\[
 (z)\bigl((w)\Theta_V\bigr)
 =\sum_{j=1}^5 (z)(a_i\widehat f_j)\otimes(v_j)a_j^V
 =\sum_{j=1}^5 (z)f_j\otimes(v_j)a_j^V.
\]
Consequently,
\[
 \bigl((w)\Theta_V\bigr)\beta_{V^{\oplus2}}^{-1}
 =\sum_{j=1}^5 f_j\otimes(v_j)a_j^V.
\]
Since $f_1,\ldots,f_5$ form a $k$-basis of $F_5$ and each $a_j^V$
is injective, $\Theta_V$ is a monomorphism of $k$-modules.

We show that
\begin{equation}\label{eq:four-ann}
 \Image(j_V)=\{z\in E(V)\mid az=0\text{ for all }a\in J^s\}.
\end{equation}
Write $K$ for the right-hand side. Since $J^sL_V=0$, we have
$\Ker(q_V)=\Image(j_V)\subseteq K$. Suppose that $(K)q_V\ne0$. Since $\End_A(S_i)=k1_{S_i}$, every
simple submodule of $S_i\otimes_k V^{\oplus5}$ has the form
$\Image(p_w)$ for some $w\ne0$. Thus $\Image(p_w)\subseteq(K)q_V$ for some
$w\ne0$. As $\Ker(q_V)\subseteq K$, the full preimage
$q_V^{-1}(\Image(p_w))$ lies in $K$. Projection onto the second
coordinate gives an isomorphism $E_w\to q_V^{-1}(\Image(p_w))$ since
$p_w$ is injective. Thus $J^sE_w=0$. Proposition~\ref{prop:comparison}
gives $(w)\Theta_V=0$, a contradiction. Hence $K=\Image(j_V)$.

We now prove that $E(V)\cong E(W)$ implies $V\cong W$. Let $V,W$
be finite-dimensional $S$-modules, and let $h\colon E(V)\to E(W)$
be an $A$-module homomorphism. By \eqref{eq:four-ann},
$(\Image(j_V))h\subseteq\Image(j_W)$. Thus the restriction of $h$ gives
a unique $B$-module homomorphism $H\colon L_V\to L_W$ with
$j_Vh=Hj_W$. On the quotients, $h$ gives a map
$1_{S_i}\otimes u$, where $u\colon V^{\oplus5}\to W^{\oplus5}$
is uniquely determined since $\End_A(S_i)=k1_{S_i}$. Hence
\[
 j_Vh=Hj_W,\qquad hq_W=q_V(1_{S_i}\otimes u).
\]
The triple $(H,h,1_{S_i}\otimes u)$ is a morphism from $\xi_V$ to
$\xi_W$. Since $p_w(1_{S_i}\otimes u)=p_{(w)u}$, the map
$E_w\to E_{(w)u}$, $(c,e)\mapsto(c,(e)h)$, gives a morphism of the
pullback extensions with end maps $H$ and $1_{S_i}$. Therefore
$([p_w^*\xi_V])\Ext_A^1(S_i,H)=[p_{(w)u}^*\xi_W]$.
Naturality of $\rho$ gives
\begin{equation}\label{eq:intertwining}
 \Theta_V F_{i,s}(H)=u\Theta_W.
\end{equation}

Choose bases of $V$ and $W$ over $k$. In these bases the entries
of $H$ lie in $R$. Replacing each entry $r$ by $\bar r$ gives a
homomorphism of $k$-modules $g\colon V^{\oplus2}\to W^{\oplus2}$. The residue map
$R\to k$ fixes scalars. Thus a change of bases gives the corresponding
change in the matrix of $g$, so $g$ is independent of the chosen bases.
Matrix reduction respects sums, identities and composition.

For $1\leq j,\ell\leq5$, write $u_{j\ell}=\iota_j u\pi_\ell$, where
$\iota_j\colon V\to V^{\oplus5}$ is the $j$th inclusion and
$\pi_\ell\colon W^{\oplus5}\to W$ is the $\ell$th projection. Apply
\eqref{eq:intertwining} to $(v)\iota_j$ and then apply
$\beta_{W^{\oplus2}}^{-1}$. The right-hand side becomes
$\sum_{\ell=1}^5 f_\ell\otimes((v)u_{j\ell})a_\ell^W$.
On the left, the entries of $H$ act on $f_j$ by postcomposition.
The equality $f_jr+F_{j-1}=f_j(\bar r\,1_Y)+F_{j-1}$ shows that
the left-hand side lies in $F_j\otimes_k W^{\oplus2}$ and maps to
$(f_j+F_{j-1})\otimes((v)a_j^V)g$ in
$(F_j/F_{j-1})\otimes_k W^{\oplus2}$. Comparing the $f_j$-coordinates
in the basis $f_1,\ldots,f_5$ of $F_5$ gives
\[
 a_j^V g=u_{jj}a_j^W\qquad(1\leq j\leq5).
\]
The first two equalities give $g=u_{11}\oplus u_{22}$, and the third
gives $u_{11}=u_{22}=u_{33}$. Put $b=u_{33}$. Comparing first
coordinates for $j=4,5$ gives $u_{44}=u_{55}=b$. The second coordinates
then give $(xv)b=x((v)b)$ and $(yv)b=y((v)b)$ for every $v\in V$.
Thus $b\colon V\to W$ is an $S$-module homomorphism.

The assignment $h\mapsto b$ respects sums, identities and composition
since $g=b\oplus b$. For an $S$-module homomorphism $b\colon V\to W$,
naturality of $j_V$ gives
$j_V E(b)=(1_Y\otimes b^{\oplus2})j_W$. Thus the assignment sends
$E(b)$ to $b$. In particular, an isomorphism $E(V)\to E(W)$ and its
inverse give inverse $S$-module homomorphisms between $V$ and $W$.

It remains to prove that $E$ preserves indecomposability. Let $V$ be
indecomposable, and take $h\in\End_A(E(V))$ with $h^2=h$.
Then $b^2=b$, so $b=0$ or $b=1_V$. If $b=0$, the entries of $H$
lie in $J(R)$. Since $J(R)$ is nilpotent and $H^2=H$, we have $H=0$
and $j_Vh=0$. The map $h$ acts as the identity on $\Image(h)$, so
$\Image(h)\cap\Image(j_V)=0$. Every simple submodule of $E(V)$ lies in
$\Image(j_V)$ by \eqref{eq:four-ann}. If $\Image(h)\ne0$, then any simple submodule of $\Image(h)$ lies in
$\Image(h)\cap\Image(j_V)=0$. Hence $h=0$. If $b=1_V$, apply the case $b=0$ to $1_{E(V)}-h$ to get
$h=1_{E(V)}$. The module $E(V)$ is nonzero by \eqref{eq:four-extension}
and has only trivial idempotent endomorphisms. Hence $E(V)$ is
indecomposable, and $E$ is a representation embedding.

Suppose now that $A$ is $\tau$-tame. If $\dim_k F\geq5$, the bimodule
$T$ constructed above has right rank $r=2\dim_k Y+5\dim_k S_i$.
By Theorem~\ref{thm:main}, for some $1\leq d\leq r^2$ there would be
infinitely many pairwise nonisomorphic indecomposable $A$-modules $X$
of dimension $d$ with $X\not\cong\tau_A(X)$. This contradicts $\tau$-tameness.
Hence $\dim_k F\leq4$, so there is an
epimorphism of $k$-modules $k^4\twoheadrightarrow F$. The isomorphism
$\alpha_Y$ in Proposition~\ref{prop:comparison} gives
$\dim_k C_{S_i}(Y)=\dim_k F$.
\end{proof}

\section{Representation dimension}\label{sec:dimensions}

Throughout this section, all $A$-modules are finite-dimensional.
Put $J=J(A)$. For a module $M$, let $\add(M)$, $\Fac(M)$ and
$\Sub(M)$ denote the full subcategories of modules isomorphic to
direct summands, quotients and submodules, respectively, of finite
direct sums of copies of $M$. We use the bound in Lemma~\ref{lem:four}
to prove that $\Fac({}_A J)$ is of finite representation type when
$A$ is tame. From the additive generators of $\Fac({}_A J)$ and
$\Sub(D(J_A))$ we construct the approximations needed for Theorem~\ref{thm:B}.

\subsection{Quotients of the radical}\label{sec:finite}

For a subcategory closed under finite direct sums and quotients,
the last nonzero radical layers determine whether the subcategory
is of finite representation type.

\begin{lemma}\label{lem:last-radical-layer}
Over an arbitrary field $k$, let $\mathcal C$ be a full subcategory of $\modu{A}$ closed under
isomorphisms, finite direct sums and quotients. Then $\mathcal C$ is of finite representation type if and only if,
for each $s\geq1$, there is a module $H_s\in\modu{A}$ such that
$J^{s-1}M$ is a quotient of $H_s$ for every indecomposable
$M\in\mathcal C$ of Loewy length $s$.
\end{lemma}

\begin{proof}
Suppose that $\mathcal C$ is of finite representation type. For each $s$, choose one
representative $M$ from each isomorphism class of indecomposable
objects of Loewy length $s$. Take $H_s$ to be the direct sum of
the modules $J^{s-1}M$ over these representatives. The
projections onto the summands give the required epimorphisms.

Conversely, suppose that the modules $H_s$ exist and that $\mathcal C$
is of infinite representation type. The duality $D$ takes
epimorphisms to monomorphisms and preserves indecomposability. Hence
the right $A$-modules isomorphic to $D(X)$ for $X\in\mathcal C$ form
a subcategory of infinite representation type closed under finite
direct sums and submodules. By \cite[Theorems~1 and~2(b)]{RingelMin},
this subcategory contains a minimal infinite submodule-closed
subcategory in which each module embeds in an indecomposable module.
Applying $D$ again gives a nonzero subcategory $\mathcal B$ of
$\mathcal C$, closed under finite direct sums. Every object of
$\mathcal B$ is a quotient of an indecomposable object of $\mathcal B$.

Since $J$ is nilpotent, we can choose a nonzero
$N\in\mathcal B$ of maximal Loewy length $s$. For each $m\geq1$, the
module $N^{\oplus m}$ belongs to $\mathcal B$. Thus there is an
epimorphism $f_m\colon M_m\twoheadrightarrow N^{\oplus m}$ with $M_m$
indecomposable in $\mathcal B$. Maximality gives $J^sM_m=0$.
Since $f_m$ is surjective,
\[
 (J^{s-1}M_m)f_m
 =J^{s-1}(N^{\oplus m})
 =(J^{s-1}N)^{\oplus m}\ne0.
\]
Hence $M_m$ has Loewy length $s$, and we have epimorphisms
\[
 H_s\twoheadrightarrow J^{s-1}M_m
 \xrightarrow{\ f_m|_{J^{s-1}M_m}\ }
 (J^{s-1}N)^{\oplus m}.
\]
It follows that $\dim_k H_s\geq m\dim_k J^{s-1}N$ for every $m\geq1$.
Since $J^{s-1}N\ne0$, taking $m=\dim_k H_s+1$ gives a contradiction.
\end{proof}

\noindent\textbf{For the next proposition and the finiteness assertion in
Remark~\ref{rem:radical-socle},} assume that $A$ is $\tau$-tame.

\begin{proposition}\label{prop:fac}
The category $\Fac({}_A J)$ is of finite representation type.
\end{proposition}

\begin{proof}
The category $\Fac({}_A J)$ is closed under isomorphisms, finite direct
sums and quotients. Let $M\in\Fac({}_A J)$ be indecomposable of Loewy
length $s$. If $s=1$, then $M$ is simple and hence a quotient of $A/J$.

Suppose $s\geq2$. Since $J^sM=0$, we may regard $M$ as an
$A/J^s$-module. Choose orthogonal primitive idempotents
$e_1,\ldots,e_n$ with sum $1$ and an epimorphism
$p\colon J^{\oplus m}\twoheadrightarrow M$. For $1\leq b\leq m$
and $1\leq i\leq n$, define
\[
 p_{i,b}\colon Je_i\longrightarrow M,\qquad
 (z)p_{i,b}=(0,\ldots,0,z,0,\ldots,0)p,
\]
where the $b$th entry is $z\in Je_i\subseteq J$ and all other entries are zero.
Let $a_i\colon J^se_i\hookrightarrow Je_i$ be the inclusion in
\eqref{eq:restriction}. Then $a_i p_{i,b}\in F_{i,s}(M)$. Since
$J^s=\bigoplus_{i=1}^n J^se_i$, we have
\[
 J^{s-1}M=((J^s)^{\oplus m})p
 =\sum_{b=1}^m\sum_{i=1}^n\Image(a_i p_{i,b}).
\]

By Lemma~\ref{lem:four}, for each $i$ we may choose four elements
$g_{i,1},\ldots,g_{i,4}$ spanning $F_{i,s}(M)$ over $k$.
By \eqref{eq:radical-layer-map}, each $g_{i,t}$ has image in
$J^{s-1}M$, and the composite
$J^{s+1}e_i\hookrightarrow J^se_i\xrightarrow{g_{i,t}}M$ is zero.
We therefore have an $A$-module homomorphism
\begin{equation}\label{eq:last-layer}
 (J^s/J^{s+1})^{\oplus4}\longrightarrow J^{s-1}M,\qquad
 (x_t+J^{s+1})_{t=1}^4\longmapsto
 \sum_{t=1}^4\sum_{i=1}^n(x_t e_i)g_{i,t}.
\end{equation}
For $x_t\in J^se_i$ and all other coordinates zero, the value of
\eqref{eq:last-layer} is $(x_t)g_{i,t}$. Thus the image of
\eqref{eq:last-layer} contains $\Image(g_{i,t})$ for each $i,t$.
Since the $g_{i,t}$ span $F_{i,s}(M)$, the image also contains every
$\Image(a_i p_{i,b})$. The sum over $i$ and $b$ is $J^{s-1}M$, so \eqref{eq:last-layer}
is surjective. In particular,
$\dim_k J^{s-1}M\leq4\dim_k(J^s/J^{s+1})$.
Lemma~\ref{lem:last-radical-layer} applies with $H_1=A/J$ and
$H_s=(J^s/J^{s+1})^{\oplus4}$ for $s\geq2$.
\end{proof}

\begin{remark}\label{rem:radical-socle}
Over any field $k$, for every $A$-module $M$ we have
\[
\begin{aligned}
 M\in\Fac({}_A J)
 &\quad\Longleftrightarrow\quad
 M\cong JX \text{ for some $A$-module $X$},\\
 M\in\Sub(D(J_A))
 &\quad\Longleftrightarrow\quad
 M\cong Z/\soc(Z) \text{ for some $A$-module $Z$}.
\end{aligned}
\]
The functor $D$ restricts to a duality between $\Fac(J_A)$ and
$\Sub(D(J_A))$. Under the hypothesis of Proposition~\ref{prop:fac},
these two categories are of finite representation type.
\end{remark}

\begin{proof}
Given an $A$-module $X$, choose an epimorphism
$p\colon A^{\oplus m}\twoheadrightarrow X$. Its restriction to
$J^{\oplus m}$ maps onto $JX$ because
$(J^{\oplus m})p=(JA^{\oplus m})p=JX$. Hence $JX\in\Fac({}_A J)$.

Conversely, for $M\in\Fac({}_A J)$ choose an epimorphism
$p\colon J^{\oplus m}\twoheadrightarrow M$. Put $N=\Ker(p)$. Since
$N\subseteq J^{\oplus m}\subseteq A^{\oplus m}$, we can take
$X=A^{\oplus m}/N$. Then
\[
 JX=\frac{JA^{\oplus m}+N}{N}=\frac{J^{\oplus m}}{N},
\]
and $x+N\mapsto(x)p$ is an $A$-module isomorphism $JX\to M$.

For an epimorphism $p\colon J_A^{\oplus m}\twoheadrightarrow Y$ of
right $A$-modules, the homomorphism
\[
 D(p)\colon D(Y)\hookrightarrow D(J_A^{\oplus m}),\qquad
 \lambda\longmapsto p\lambda,
\]
is injective. Restricting to each summand gives
$D(J_A^{\oplus m})\cong D(J_A)^{\oplus m}$, so
$D(Y)\in\Sub(D(J_A))$. Conversely, a monomorphism
$i\colon M\hookrightarrow D(J_A)^{\oplus m}$ gives the epimorphism
\[
 J_A^{\oplus m}\xrightarrow{\sim}D\bigl(D(J_A)^{\oplus m}\bigr)
 \xrightarrow{D(i)}D(M).
\]
The first arrow sends $(z_1,\ldots,z_m)$ to the functional
$(\lambda_1,\ldots,\lambda_m)\mapsto\sum_{j=1}^m(z_j)\lambda_j$.
Consequently $D$ restricts to the claimed duality.

By the first equivalence applied to right $A$-modules, the objects of
$\Fac(J_A)$ are the modules isomorphic to $YJ$ for some right $A$-module
$Y$. The objects of $\Sub(D(J_A))$ are therefore the modules isomorphic
to $D(YJ)$ for such $Y$. Applying $D$ to the inclusion $YJ\hookrightarrow Y$
gives the epimorphism
\begin{equation}\label{eq:dual-radical}
 D(Y)\twoheadrightarrow D(YJ),\qquad
 \lambda\longmapsto\lambda|_{YJ}.
\end{equation}
For $a\in J$ and $y\in Y$, we have $(y)(a\lambda)=(ya)\lambda$.
Thus $\lambda|_{YJ}=0$ if and only if $J\lambda=0$. The kernel of
\eqref{eq:dual-radical} is $\soc(D(Y))$, and restriction gives the
$A$-module isomorphism
\[
 D(Y)/\soc(D(Y))\xrightarrow{\sim}D(YJ),\qquad
 \lambda+\soc(D(Y))\longmapsto\lambda|_{YJ}.
\]
Every $A$-module $Z$
is naturally isomorphic to $D(D(Z))$, so the second equivalence follows.

Under the hypothesis of Proposition~\ref{prop:fac}, $A\opp$ is
$\tau$-tame by Corollary~\ref{cor:descent}. Applying the proposition
to $A\opp$ shows that $\Fac(J_A)$ is of finite representation type.
The assertion for $\Sub(D(J_A))$ follows by duality.
\end{proof}

\subsection{Right approximations}\label{sec:approximations}

The following construction works over any field $k$. It uses additive
generators of $\Fac({}_A J)$ and $\Sub(D(J_A))$ to give right
approximations.

For any $A$-module $G$ and $G_0\in\add(G)$, a homomorphism
$h\colon G_0\to X$ is a right $\add(G)$-approximation if and only if
\[
 \Hom_A(G,h)\colon\Hom_A(G,G_0)\longrightarrow\Hom_A(G,X),
 \qquad f\longmapsto fh,
\]
is surjective. Indeed, if every homomorphism $G\to X$ factors through
$h$, the same holds for finite direct sums of copies of $G$ and their
direct summands.

\begin{proposition}\label{prop:approx}
Let $V$ and $U$ be $A$-modules such that
$\Fac({}_A J)=\add(V)$ and $\Sub(D(J_A))=\add(U)$. Put
$G={}_A A\oplus D(A_A)\oplus V\oplus U$.
For every $X\in\modu{A}$, there is a short exact sequence
\begin{equation}\label{eq:approx}
 0\longrightarrow G_1\xrightarrow{d}G_0\xrightarrow{h}X
 \longrightarrow0
\end{equation}
with $G_0,G_1\in\add(G)$ such that applying $\Hom_A(G,-)$ gives an
exact sequence.
\end{proposition}

\begin{proof}
Finite direct sums of the required sequences reduce the proof to
indecomposable $X$.
If $X$ is injective, take $G_0=X$, $h=1_X$ and $G_1=0$.
Assume that $X$ is noninjective. Let $L$ be the trace of $J$ in $X$, that is, the image of
\[
 J\otimes_k\Hom_A(J,X)\longrightarrow X,\qquad z\otimes f\longmapsto(z)f.
\]
Then $L\in\Fac({}_A J)$. Given $W\in\Fac({}_A J)$ and $f\colon W\to X$,
choose an epimorphism $J^{\oplus m}\twoheadrightarrow W$. The
restriction of the composite $J^{\oplus m}\to W\xrightarrow{f}X$
to each summand has image in $L$. Hence $\Image(f)\subseteq L$.
An epimorphism $A^{\oplus m}\twoheadrightarrow X$ maps $J^{\oplus m}$
onto $JX$, so $JX\subseteq L$. Thus $X/L$ is semisimple.
Let $q\colon X\twoheadrightarrow X/L$ be the quotient map.

Put $U_X=U\otimes_k\Hom_A(U,X)$ and define $u\colon U_X\to X$ by
$(w\otimes f)u=(w)f$. Since $U_X\in\add(U)$ and each $f\colon U\to X$
factors through $u$ by $w\mapsto w\otimes f$, the map $u$ is a right
$\add(U)$-approximation. Choose a complement $T$ of $\Image(uq)$ in $X/L$.
Write $i\colon T\hookrightarrow X/L$ for the inclusion. Take a projective
cover $v\colon P\twoheadrightarrow T$ and lift $vi$ to $p\colon P\to X$,
so that $pq=vi$. Since $T$ is semisimple and $v$ is a projective cover,
\[
 X/L=\Image(uq)\oplus\Image(pq),\qquad \Ker(pq)=\Ker(v)=JP.
\]

Define
\begin{equation}\label{eq:h}
 h\colon L\oplus U_X\oplus P\longrightarrow X,\qquad
 (l,y,z)\longmapsto l+(y)u+(z)p.
\end{equation}
The map $h$ is surjective because $L\subseteq\Image(h)$ and $(\Image(h))q=X/L$.
Every homomorphism $V\to X$ factors through the inclusion of $L$.
Every homomorphism $U\to X$ factors through $u$, and every homomorphism
$A\to X$ factors through $h$ by projectivity.

It remains to consider a homomorphism $f\colon D(A_A)\to X$. Fix an
indecomposable injective summand $I$ of $D(A_A)$, and let $S$ be its socle.
The restriction of $f$ to $I$ cannot be a split monomorphism, since
$X$ is indecomposable and noninjective. The quotient $I\to I/S$ is
left almost split \cite[Chapter~V, Corollary~1.6(b)]{ARS}, so this
restriction factors through $I/S$. Hence the composite
$S\hookrightarrow I\hookrightarrow D(A_A)\xrightarrow{f}X$ is zero.
Taking the sum over $I$ shows that $\soc(D(A_A))\subseteq\Ker(f)$.
By \eqref{eq:dual-radical} with $Y=A_A$, this socle is the kernel
of the restriction $r\colon D(A_A)\twoheadrightarrow D(J_A)$.
There is therefore a unique $f'\colon D(J_A)\to X$ such that $rf'=f$.
Since $D(J_A)\in\add(U)$, the map $f'$ factors through $u$.
This proves that $h$ is a right $\add(G)$-approximation.

For $(l,y,z)\in\Ker(h)$, applying $q$ to \eqref{eq:h} gives
$(y)uq+(z)pq=0$. Since $\Image(uq)\cap\Image(pq)=0$, we have $y\in\Ker(uq)$
and $z\in JP$. Projection onto the last two summands gives the isomorphism
\[
 \Ker(h)\xrightarrow{\sim}\Ker(uq)\oplus JP,\qquad
 (l,y,z)\longmapsto(y,z),
\]
with inverse $(y,z)\mapsto(-(y)u-(z)p,y,z)$.
As a submodule of $U_X$, the module $\Ker(uq)$ belongs to
$\Sub(D(J_A))=\add(U)$. Since $P$ is a direct summand of $A^{\oplus m}$
for some $m$, we have $JP\in\add({}_A J)\subseteq\add(V)$.
Thus $G_1=\Ker(h)$ and $G_0=L\oplus U_X\oplus P$ both belong to $\add(G)$.
Take $d$ in \eqref{eq:approx} to be the inclusion. The functor
$\Hom_A(G,-)$ is left exact, and $\Hom_A(G,h)$ is surjective.
Thus applying $\Hom_A(G,-)$ to \eqref{eq:approx} gives an exact sequence.
\end{proof}

\noindent\textbf{For the rest of this section, assume that $A$ is of tame
representation type.} We now apply the approximation criterion to the
two additive generators obtained in Subsection~\ref{sec:finite}.

\begin{proposition}\label{prop:generator}
The categories
$\Fac({}_A J)$ and $\Sub(D(J_A))$ are of finite representation type. For any $A$-modules
$V$ and $U$ satisfying
\[
 \add(V)=\Fac({}_A J),\qquad \add(U)=\Sub(D(J_A)),
\]
the module
\[
 G={}_A A\oplus D(A_A)\oplus V\oplus U
\]
is a generator-cogenerator with $\gldim\End_A(G)\leq3$. If $A$ is of
infinite representation type, then $\repdim A=3$.
\end{proposition}

\begin{proof}[Proof of Proposition~\ref{prop:generator} and Theorem~\ref{thm:B}]
By Theorem~\ref{thm:main}, $A$ is $\tau$-tame. By Proposition~\ref{prop:fac}
and Remark~\ref{rem:radical-socle}, $\Fac({}_A J)$ and $\Sub(D(J_A))$
are of finite representation type. Thus each category has an additive
generator. The summands ${}_A A$ and $D(A_A)$ make $G$ a
generator-cogenerator. By Proposition~\ref{prop:approx}, every
$A$-module has a short exact sequence of the form \eqref{eq:approx}
which remains exact under $\Hom_A(G,-)$. The approximation criterion
\cite[Lemma~2.1]{EHIS} gives $\gldim\End_A(G)\leq3$. Since $G$ is
a generator-cogenerator, $\repdim A\leq\gldim\End_A(G)\leq3$,
as asserted in Theorem~\ref{thm:B}. If $A$ is of infinite
representation type, Auslander's characterization \cite{Auslander}
gives $\repdim A>2$. Thus $\repdim A=3$.
\end{proof}

\subsection{Homological applications}\label{sec:applications}

Choose $V$ and $U$ as in Proposition~\ref{prop:generator}, and put
$H=V\oplus A/J$. For every $A$-module $M$, the radical sequence
$0\to JM\to M\to M/JM\to0$ has both end terms in $\add(H)$.
Remark~\ref{rem:radical-socle} gives the dual statement for the socle
sequence with $U\oplus A/J$. Thus $\modu A$ has extension dimension at most one. By definition,
the extension dimension is zero precisely when $\modu A=\add(N)$
for some $A$-module $N$, that is, when $A$ is of finite representation type.

Write $\mathrm D_{\mathrm{sg}}(A)$ for the singularity category, the
quotient of $\mathrm D^b(\modu A)$ by the bounded homotopy category of
finitely generated projective modules. We use Rouquier dimension for
triangulated categories \cite[Definition~3.6]{Rouquier}.

\begin{corollary}\label{cor:derived-dimension}
The bounded derived category and the singularity category satisfy
\[
 \dim\mathrm D^b(\modu A)\leq2,\qquad
 \dim\mathrm D_{\mathrm{sg}}(A)\leq1.
\]
Both bounds are sharp.
\end{corollary}

\begin{proof}
For a bounded complex $X^\bullet$, let $B^\bullet$ be the subcomplex
of boundaries of $JX^\bullet$. The filtration
$0\subseteq B^\bullet\subseteq JX^\bullet\subseteq X^\bullet$
has three factors. The first two have zero differentials and terms in
$\add(V)=\Fac({}_A J)$, since this category is closed under quotients.
The last factor splits over $A/J$ into shifts of its cohomology modules
in the derived category. Thus $\mathrm D^b(\modu A)=\langle H\rangle_3$.
Every object in the singularity category is isomorphic to a shift of the
image of a module \cite[proof of Proposition~3.7]{Rouquier}, so the radical
sequence gives $\mathrm D_{\mathrm{sg}}(A)=\langle q(H)\rangle_2$ for
the quotient functor $q$. Both bounds are attained by the exterior
algebra on two generators \cite[Theorem~4.1]{Rouquier}, which is
special biserial and hence tame \cite[Introduction]{EHIS}.
\end{proof}

For self-injective $A$, the stable module category is equivalent to
$\mathrm D_{\mathrm{sg}}(A)$ as a triangulated category \cite[Section~3.3]{Rouquier}.
Corollary~\ref{cor:derived-dimension} therefore gives the bound in
Yoshiwaki's conjecture \cite[p.~107, Conjecture~2]{Yoshiwaki}.
For connected self-injective $A$ of infinite representation type,
the stable dimension is exactly one by Proposition~\ref{prop:generator}
and \cite[Corollary~3.11]{Yoshiwaki}.

We next consider the finitistic dimension conjecture.

\begin{corollary}\label{cor:findim}
Both $A$ and $A\opp$ satisfy the finitistic dimension conjecture.
\end{corollary}

\begin{proof}
Since $A$ and $A\opp$ are tame, both have representation dimension
at most three by Theorem~\ref{thm:B}. The assertion follows from
the theorem of Igusa and Todorov \cite{IT}.
\end{proof}

Zhang and Zhou \cite[Proposition~4.1]{ZZARC} proved the
Auslander--Reiten conjecture for algebras of representation dimension
at most three. Theorem~\ref{thm:B} gives the following consequence.

\begin{corollary}\label{cor:projectivity}
$A$ satisfies the Auslander--Reiten conjecture.
\end{corollary}

By Corollary~\ref{cor:findim}, both $A$ and $A\opp$ have finite
little finitistic dimension. Hence $A$ satisfies the Gorenstein
symmetry conjecture by \cite[Proposition~6.10]{AR1991}.

\begin{corollary}\label{cor:gorenstein-symmetry}
$A$ satisfies the Gorenstein symmetry conjecture.
\end{corollary}

\paragraph*{Acknowledgements.}The authors used ChatGPT to help organize parts of the manuscript and revise the English exposition and mathematical notation.
All four authors reviewed the material they adopted, made the final revisions, and take responsibility for the entire article.

\Needspace{12\baselineskip}

\bigskip
\noindent\begin{minipage}{\textwidth}
\normalsize
\setlength{\parindent}{1.5em}
Liang Chen: School of Mathematical Sciences, Capital Normal
University, 100048 Beijing, China\par
\textit{E-mail address:}
\href{mailto:2210501003@cnu.edu.cn}{\texttt{2210501003@cnu.edu.cn}}\par
\end{minipage}

\bigskip
\noindent\begin{minipage}{\textwidth}
\normalsize
\setlength{\parindent}{1.5em}
Zhenxian Chen: School of Mathematical Sciences, Capital Normal
University, 100048 Beijing, China\par
\textit{E-mail address:}
\href{mailto:czx18366459216@163.com}{\texttt{czx18366459216@163.com}}\par
\end{minipage}

\bigskip
\noindent\begin{minipage}{\textwidth}
\normalsize
\setlength{\parindent}{1.5em}
Yongyi Wen: School of Mathematical Sciences, Capital Normal
University, 100048 Beijing, China\par
\textit{E-mail address:}
\href{mailto:2240501002@cnu.edu.cn}{\texttt{2240501002@cnu.edu.cn}}\par
\end{minipage}

\bigskip
\noindent\begin{minipage}{\textwidth}
\normalsize
\setlength{\parindent}{1.5em}
Weiheng Xia: School of Mathematical Sciences, Capital Normal
University, 100048 Beijing, China\par
\textit{E-mail address:}
\href{mailto:xia.weiheng@qq.com}{\texttt{xia.weiheng@qq.com}}
\end{minipage}

\end{document}